\documentclass[10pt,a4paper,leqno]{amsart}

\usepackage{amsmath,amssymb,amsthm,euscript}

\def\eps
{\varepsilon}

\def\R{{\mathbb R}}
\def\N{{\mathbb N}}

\def\virgp{\raise 2pt\hbox{,}}

\def\bu{{\bf u}}
\def\({\left(}
\def\){\right)}

\def\Eq#1#2{\mathop{\sim}\limits_{#1\rightarrow#2}}
\def\Tend#1#2{\mathop{\longrightarrow}\limits_{#1\rightarrow#2}}

\def\d{{\partial}}
\def\e{\varepsilon}

\def\U{{\mathcal U}^\e}
\def\V{{\tt V}^\e}

\def\u{{\tt u}}

\def\F{\mathcal F}
\def\M{{\mathcal M}^\e}
\def\D{{\tt D}^\e}
\def\O{\mathcal O}

\theoremstyle{plain}
\newtheorem{theo}{Theorem}[section]
\newtheorem{lem}[theo]{Lemma}
\newtheorem{cor}[theo]{Corollary}
\newtheorem{prop}[theo]{Proposition}

\theoremstyle{definition}
\newtheorem*{defin}{Definition}
\newtheorem*{nota}{Notation}

\theoremstyle{remark}

\newtheorem*{rema*}{Remark}

\numberwithin{equation}{section}

\begin{document}

\title[Semiclassical Hartree equation
with harmonic potential]{(Semi)classical limit of the Hartree equation
with harmonic potential}
\author[R. Carles]{R{\'e}mi Carles}
\address[R. Carles\footnote{On leave from MAB, Universit\'e Bordeaux
1.}]{IRMAR, Universit\'e de Rennes 1\\ Campus de
Beaulieu\\ 35~042 Rennes cedex\\ France}
\email{remi.carles@math.univ-rennes1.fr}
\author[N. Mauser]{Norbert J. Mauser}
\address[N.J.  Mauser]{Wolfgang Pauli Institute c/o Inst. f. Math. \\
Universit\"at Wien\\
Nordbergstr. 15 \\ A 1090 Wien\\ Austria}
\email{mauser@courant.nyu.edu}
\author[H. P. Stimming]{Hans Peter Stimming}
\address[H. P. Stimming]{Wolfgang Pauli Institute, Wien and ENS Lyon, France}
\email{hans.peter.stimming@univie.ac.at}
\thanks{The authors acknowledge support by the Austrian START award project
(FWF, contract No. Y-137-TEC) of N.J.M. and
by the Wissenschaftskolleg (doctoral school) ``Differential Equations''
(FWF, contract No. W8)
as well as the European network HYKE, funded  by the EC as contract
HPRN-CT-2002-00282}
\begin{abstract}
Nonlinear Schr\"odinger Equations (NLS) of the Hartree type occur
in the modeling of quantum semiconductor devices. Their
"semiclassical" limit of vanishing (scaled) Planck constant is
both a mathematical challenge and practically relevant when
coupling quantum models to classical models.
 With the aim of describing the semi-classical limit of the 3D
Schr\"o\-dinger--Poisson system with an additional harmonic
potential, we study some semi-classical limits of the Hartree
equation with harmonic potential in space dimension $n\geq~2$. The
harmonic potential is confining, and causes focusing periodically
in time. We prove asymptotics in several cases, showing different
possible nonlinear phenomena according to the interplay of the
size of the initial data and the power of the Hartree potential.
In the case of the 3D Schr\"odinger--Poisson system with harmonic
potential, we can only give a formal computation since the need of
modified scattering operators for this long range scattering case
goes beyond current theory. \\
We also deal with the case of an additional "local" nonlinearity
given by a power of the local density - a model that is relevant
when incorporating the Pauli principle in the simplest model given
by the "Schr\"odinger-Poisson-X$\alpha$ equation". Further we
discuss the connection of our WKB based analysis to the Wigner
function approach to semiclassical limits.
\end{abstract}
\maketitle

\section{Introduction}
\label{sec:intro}

Nonlinear Schr\"odinger Equations (NLS) are important both for
many different applications as well as a source of rich
mathematical theory with several hard challenges still open. The
NLS in the most common meaning contains a "local" nonlinearity
given by a power of the local density, in particular the
(de)focusing "cubic" NLS which arises e.g. in nonlinear optics or
for Bose Einstein condensates. In 1-d this NLS is an integrable
system and the "semi-classical limit" ("high wave number limit")
can be performed by methods of inverse scattering (see e.g.
\cite{JLM} and \cite{KMM} for results on the defocusing and
focusing case). A class of NLS with a "nonlocal" nonlinearity that
we call "Hartree type" occur in the modeling of quantum
semiconductor devices. Their "semi-classical" limit of vanishing
(scaled) Planck constant is both a mathematical challenge and
practically relevant when coupling quantum models to
classical models.\\
Incorporating the Pauli principle for fermions in the simplest
possible model yields the case of a Hartree equation with an
additional "local" nonlinearity given by a power of the local
density, the "Schr\"odinger-Poisson-X$\alpha$ equation" (see
\cite{MauserAML}).\\

In this paper we deal with the ``semi-classical limit" of nonlinear
Schr\"odinger equations of  Hartree type, with a harmonic
potential and a ``weak" nonlinearity which is a convolution of the density with
a more or less singular potential.

In three space dimensions, for the case where we convolute
with the Newtonian potential $1/|x|$,
the Hartree equation is the Schr\"odinger--Poisson
system with  harmonic potential :
\begin{equation}\label{eq:s-p}
\left\{
\begin{aligned}
i\e \d_t \u^\e +\frac{1}{2}\e^2\Delta \u^\e &= \frac{|x|^2}{2}\u^\e +
V(x)\u^\e \, ,\\
 \Delta V &= |\u^\e|^2\, ,\\
\u^\e_{\mid t=0}& =
\u^\e_0\, ,
\end{aligned}
\right.
\end{equation}
with $x\in \R^3$.

This equation arises typically if we consider the quantum mechanical
time evolution of electrons in the mean field approximation of the
many body effects, modeled by the Poisson equation, with
a confinement modeled by the quadratic potential of the harmonic oscillator.

The limit $\e\to 0$ in such a quantum model corresponds to a ``classical
limit" of vanishing Planck constant $\hbar = \e \to 0$. We adopt the terminology
``semi-classical limit'' for what should properly be called ``classical limit''
(see the discussion in \cite{ZZM}), the term ``semi-classical'' being actually more appropriate
for the situation of the homogenization limit from a Schr\"odinger equation with
periodic potential (see e.g. \cite{BMP2}).

The problem of the mathematically rigorous ``classical limit" of the
Schr\"odinger-Poisson system is highly nontrivial.
First results of weak limits $\e\to 0$ to the Vlasov-Poisson system
where given in \cite{LionsPaul} and \cite{MM1} using Wigner transform
techniques for the ``mixed state case", where additional strong assumptions
on the initial data can be imposed (which are necessary to guarantee a
uniform $L^2$ bound on the Wigner function).
In \cite{ZZM} this assumption could be removed for the 1-d case
and the classical limit for the "pure state" case could be performed,
where the notorious problem of non-uniqueness of the Vlasov-Poisson system
with measure valued initial data reappears.
For an overview of this kind of ``semi-classical limits" of Hartree equations
see \cite{MauserForges}.
For an introduction to Wigner transforms and their comparison to WKB methods
for the linear case see \cite{GMMP} and \cite{SMM}.

Up to a constant,
\eqref{eq:s-p} is equivalent to the Hartree equation
\begin{equation}\label{eq:r3.1}
i\e \d_t \u^\e +\frac{1}{2}\e^2\Delta \u^\e = \frac{|x|^2}{2}\u^\e +
\(|x|^{-1}\ast |\u^\e|^2\)\u^\e \quad ; \quad \u^\e_{\mid t=0} =
\u^\e_0\, .
\end{equation}
We restrict our attention to small data cases with
$\u^\e_0 = \e^{\alpha/2} f$, where $f$ is independent of $\e$ and
$\alpha \geq 1$.

Notice that we can allow for more general data with initial plane oscillations,
\begin{equation}
\label{eq:plane}
\bu^\e_{\mid t=0}= \e^{\alpha/2} f(x) e^{i\frac{x.\xi_0}{\e}}
 \quad \mbox{for} \quad \xi_0\in\R^3,
\end{equation}
since the change of variables given in \cite{CaIHP}
\begin{equation}\label{eq:Luc}
\bu^\e(t,x)=\u^\e(t,x-\xi_0 \sin t)e^{i\left( x-
\frac{\xi_0}{2}\sin t\right).\xi_0 \cos t /\e}\, ,
\end{equation}
yields the solution of \eqref{eq:r3.1}.  This
change of variable  could also be used in Equation~\eqref{eq:r3.3}
below and hence our results also hold for the more general
$\e$-dependent class of data \eqref{eq:plane}.

Note that ``small data'' can be equivalently written as ``small nonlinearity'', since
with the change of the unknown $u^\e = \e^{-\alpha/2}\u^\e$, \eqref{eq:r3.1} becomes
\begin{equation}\label{eq:r3.2}
i\e \d_t u^\e +\frac{1}{2}\e^2\Delta u^\e = \frac{|x|^2}{2}u^\e +
\e^{\alpha}\(|x|^{-1}\ast |u^\e|^2\)u^\e \quad ; \quad u^\e_{\mid t=0} =
f\, .
\end{equation}
We will consider the more general ``semi-classical Hartree equation''
\begin{equation}\label{eq:r3.3}
i\e \d_t  u^\e +\frac{1}{2}\e^2\Delta u^\e = \frac{|x|^2}{2}u^\e +
\e^\alpha \(|x|^{-\gamma}\ast |u^\e|^2\)u^\e \quad ; \quad u^\e_{\mid t=0} =
f\, ,
\end{equation}
with $\gamma>0$, $\alpha \geq 1$ and $x\in\R^n$, where the space
dimension $n\geq 2$ may be different from $3$. \\
The first point to notice is that in the linear case, the harmonic potential causes
focusing at the origin (resp. at $(-1)^k \xi_0$ in the case
\eqref{eq:Luc}) at
times $t=\pi/2 +k\pi$, for any $k\in \N$. The solution $u^\e_{\rm
free}$ of the  linear equation
\begin{equation}\label{eq:harmo}
i\e \d_t  u^\e_{\rm free} +\frac{1}{2}\e^2\Delta u^\e_{\rm free} =
\frac{|x|^2}{2}u^\e_{\rm free}  \quad ; \quad u^\e_{\rm{free}\mid t=0} =
f\, ,
\end{equation}
is initially of size $\O(1)$. At time $t=\pi/2$, the solution focuses
at the origin and is of order $\O(\e^{-n/2})$; it is of order $\O(1)$
for $t=\pi$, and so on (for a more precise
analysis, see  \cite{CaIHP}). This phenomenon is easy to read from Mehler's
formula (see e.g. \cite{Feyn,HormanderQuad}):  for $0<t<\pi$,
we have
\begin{equation}\label{eq:mehler}
u^\e_{\rm
free}(t,x)=\frac{e^{-in\frac{\pi}{4}}}{(2\pi
\e\sin  t)^{n/2}}\int_{\R^n}e^{\frac{i}{ \e\sin t}
\left(\frac{x^2+y^2}{2}\cos t -x\cdot y \right)}f(y)dy \, .
\end{equation}
Essentially, one can apply a stationary phase formula for $t\in
]0,\pi/2[\cup ]\pi/2,\pi[$ ($u^\e_{\rm
free}$ is $\O(1)$), while it is not possible at $t=\pi/2$ ($u^\e_{\rm
free}$ is $\O(\e^{-n/2})$).
Following the same approach as in \cite{Ca2}, we get the following
distinctions:\\

\begin{center}
\begin{tabular}[c]{l|c|c}
 &$\alpha >\gamma$  & $\alpha =\gamma$  \\
\hline
$\alpha >1$ &Linear WKB, & Linear WKB,\\
&linear focus   & nonlinear focus \\
\hline
$\alpha =1$&Nonlinear WKB, &Nonlinear WKB,\\
&linear focus   & nonlinear focus \\
\end{tabular}
\end{center}
\bigskip

The expression ``linear WKB'' means that the nonlinear Hartree
interaction term is negligible away from the focus (when the WKB
approximation is valid); ``linear focus''
means that the nonlinearity is negligible near the focus; the
WKB r\'egime (resp. the focus) is ``nonlinear'' when the Hartree term
has a leading order influence away from (resp. in the neighborhood of)
the focus, in the limit $\e \to 0$. This terminology follows \cite{HK87}.

We did not obtain a rigorous description of the case
$\alpha =\gamma=1$, which corresponds to the Schr\"odinger--Poisson
system \eqref{eq:s-p} when $n=3$. This problem seems out of reach for the
methods currently available in this field.
On the other hand, we study rigorously the three other cases in an exhaustive
way:

 In Section~\ref{sec:small}, we prove that the Hartree
term has no influence at leading order when $\alpha >\gamma=1$. Back
to \eqref{eq:r3.1}, this shows that initial data of size
$\e^{\alpha/2}$ with $\alpha >1$ yield a linearizable solution. The
expected critical size is $\sqrt\e$; this heuristic is reinforced by the
next three sections.

In Section~\ref{sec:nlwkb}, we study the case $\alpha =1 >\gamma$. We
prove that the nonlinear term must be taken into account to describe
the solution $u^\e$. It is so through a slowly oscillating phase
term. On the other hand, no nonlinear effect occurs at leading order
near the focus.

In Section~\ref{sec:scattering}, we show that when $\alpha =\gamma
>1$, nonlinear effects occur at leading order at the focuses,
while they are negligible elsewhere. This phenomenon is the same as in
\cite{CaIHP} for the nonlinear Schr\"odinger equation; each focus
crossing is described in terms of the scattering operator associated
to the Hartree equation
\begin{equation}\label{eq:r3gen}
i\d_t \psi + \frac{1}{2}\Delta \psi = \( |x|^{-\gamma}\ast
|\psi|^2\)\psi\, .
\end{equation}

In Section~\ref{sec:formal}, we perform a formal computation suggested
by the results of Sections~\ref{sec:nlwkb} and
\ref{sec:scattering}. This can be seen as a further evidence that
nonlinear effects are always relevant in the case $\alpha =\gamma =1$,
along with a precise idea of the nature of these nonlinear effects,
which we expect to be true.
We add a brief discussion of the case of an additional local
nonlinearity in the equation and some remarks on the Wigner measures
in view of the ill-posedness results of \cite{CaWigner}.

This program is very similar to the one achieved in \cite{Ca2}. We
want to underscore at least two important differences. First, we have
to adapt the notion of oscillatory integral to incorporate the
presence of the harmonic potential (see
Section~\ref{sec:osc}). Second, the power-like nonlinearity treated in
\cite{Ca2} is replaced by a Hartree-type nonlinearity. This yields different and
less technical proofs (we do not use Strichartz estimates in
Sections~\ref{sec:small} and \ref{sec:nlwkb}), and makes a more
complete description of the above table possible; the case ``nonlinear WKB,
linear focus'' was treated very partially in \cite{Ca2}, due to the
lack of regularity of the map $z\mapsto |z|^{2\sigma}z$ for small
$\sigma >0$. This technical difficulty does not occur in the present
case, and the main result of Section~\ref{sec:nlwkb}
(Proposition~\ref{prop:nlwkb}) is proved with no restriction.

The content of this article is as explained above, plus a paragraph
dedicated to a quick review of the facts we will need about the
Cauchy problem \eqref{eq:r3.3} (Section~\ref{sec:cauchy}).

We will use the following notation throughout this paper.
\begin{nota}
If $(a^\e)_{\e \in]0,1]}$ and $(b^\e)_{\e \in]0,1]}$
are two families of numbers, we write
$$a^\e \lesssim b^\e$$
if there exists $C$ independent of $\e \in ]0,1]$ such that for
any $\e \in ]0,1]$, $a^\e \leq C b^\e$.
\end{nota}

\section{The Cauchy problem}
\label{sec:cauchy}
Before studying semi-classical limits, we recall some known facts
about the initial value problem \eqref{eq:r3.3}. We will always assume
that the initial datum $f$ is in the space $\Sigma$ defined by
\begin{equation*}
\Sigma := \left\{ \phi \in H^1(\R^n)\ ;\ \|\phi\|_\Sigma :=
\|\phi\|_{L^2} + \|x \phi\|_{L^2} + \|\nabla \phi\|_{L^2} <+\infty \right\}\, .
\end{equation*}
This space is natural in the case of Schr\"odinger equations with
harmonic potential, since $\Sigma$ is the domain of $\sqrt{-\Delta +|x|^2}$
(see for instance \cite{Oh}). Local existence results for
\eqref{eq:r3.3} follow for instance from Strichartz inequalities (one
can do without these inequalities, see \cite{Oh}). Global
existence results then stem from conservation laws (see
\eqref{eq:conservations} below). From
Mehler's formula \eqref{eq:mehler}, Strichartz type estimates
are available for
\begin{equation*}
e^{-i\frac{t}{2\e}(-\e^2\Delta+x^2)}=:\U(t)\, .
\end{equation*}

\begin{defin}\label{def:adm}
 Let $n\geq 2$. A pair $(q,r)$ is {\bf admissible} if $2\leq r
  <\frac{2n}{n-2}$ (resp. $2\leq r<
  \infty$ if $n=2$)
  and
$$\frac{2}{q}=\delta(r)\equiv n\left( \frac{1}{2}-\frac{1}{r}\right).$$
\end{defin}
Following \cite{CaIHP}, we have the following scaled Strichartz inequalities:
 \begin{prop}\label{prop:strichartz}
Let $I$ be a \emph{finite} time interval.\\
$(1)$ For any admissible pair $(q,r)$,
  there exists $C_r(I)$ such that
\begin{equation}\label{eq:strichlib}
   \e^{\frac{1}{q}} \left\| \U(t)\phi\right\|_{L^q(I;L^r)}\leq C_r(I)
   \|\phi\|_{L^2}\, .
  \end{equation}
$(2)$ For any admissible pairs $(q_1,r_1)$ and $(q_2,r_2)$, there
   exists $C_{r_1,r_2}(I)$ such that
\begin{equation}\label{eq:strichnl}
      \e^{\frac{1}{q_1}+\frac{1}{q_2}}\left\| \int_{I\cap\{s\leq
      t\}} \U(t-s)F(s)ds
      \right\|_{L^{q_1}(I;L^{r_1})}\leq C_{r_1,r_2}(I) \left\|
      F\right\|_{L^{q'_2}(I;L^{r'_2})}\, .
    \end{equation}
The above constants are independent of $\e$.
\end{prop}
The main result of this section follows from \cite{Caz,Ginibre}. Denote
\begin{equation*}
Y (I) =\{ \phi\in C(I,\Sigma)\ ;\  \phi, \, |x| \phi,\, \nabla_x \phi \in
L^q_{\rm loc}(I,L^r_x), \ \forall
(q,r)\textrm{ admissible}\}\, .
\end{equation*}
\begin{prop}\label{prop:cauchy}
Fix $\eps \in ]0,1]$ and let $f \in \Sigma$. Then \eqref{eq:r3.3} has
a unique solution $u^\eps \in
Y(\R)$. Moreover, the following quantities are independent of time:
\begin{equation}\label{eq:conservations}
\begin{aligned}
\text{Mass: }& \|u^\e(t)\|_{L^2}\, ,\\
\text{Energy: }& \frac{1}{2}\|\e\nabla_x
u^\e(t)\|^2_{L^2}  + \frac{1}{2}\| x u^\e(t)\|_{L^2}^2 + \e \int_{\R^n}
\(|x|^{-\gamma}\ast |u^\e|^2\)|u^\e(t,x)|^2dx\\
\end{aligned}
\end{equation}
\end{prop}
It was noticed in \cite{CaIHP} that this result can be retrieved very
simply thanks to the following lemma, which we will use to prove asymptotics.
\begin{lem}[\cite{CaIHP}]\label{lem:operators}
Define the operators
\begin{equation}\label{eq:op}
J^\e(t) = \frac{x}{\e}\sin  t
-i \cos t \nabla_x\quad ;\quad
H^\e(t)=  x\cos  t+i \e \sin  t \nabla_x\, .
\end{equation}
$J^\e$ and  $H^\e$ satisfy the following properties.\\
$\bullet$ They are Heisenberg observables:
\begin{equation}\label{eq:Heisenberg}
J^\e(t) = -i\,\U(t)\nabla_x \U(-t)\quad ; \quad H^\e(t)= \U(t)x\, \U(-t)\, .
\end{equation}
$\bullet$ The commutation relation:
\begin{equation}\label{eq:commut}
\left[J^\e(t), i \e\d_t +\frac{\e^2}{2}\Delta
-\frac{|x|^2}{2}
\right]=\left[H^\e(t), i \e\d_t +\frac{\e^2}{2}\Delta
-\frac{|x|^2}{2}
\right]=0\,.
\end{equation}
$\bullet$ Denote $M^\e(t)= e^{-i \frac{x^2}{2\e}\tan
t}$, and
$Q^\e(t) = e^{i \frac{x^2}{2\e}\cot  t}$,
then
\begin{equation}\label{eq:factor}
J^\e (t) =-i\cos t\,  M^\e(t)\nabla_x  M^\e(-t)\quad ;\quad
H^\e (t) =i\e\sin t \, Q^\e(t)\nabla_x  Q^\e(-t)\, .
\end{equation}
$\bullet$ The modified Sobolev inequalities.
Let $2\leq r\leq \frac{2n}{n-2}$ ($2\leq r <\infty$ if
$n=2$); there exists $C_r$ independent of $\e$ such
that, for any $\phi \in \Sigma$,
\begin{equation}\label{eq:sobolev}
\begin{aligned}
\|\phi\|_{L^r}  &\leq C_r |\cos t|^{-\delta(r)}
\|\phi\|_{L^2}^{1-\delta(r)}
 \| J^\e(t)\phi)\|_{L^2}^{\delta(r)}\, ,
\\
\|\phi\|_{L^r} & \leq C_r |\e\sin t|^{-\delta(r)}
\|\phi\|_{L^2}^{1-\delta(r)}
\| H^\e(t)\phi\|_{L^2}^{\delta(r)}\,.
\end{aligned}
\end{equation}
$\bullet$ Action on nonlinear Hartree term: for $\phi=\phi(t,x)$,
\begin{equation}\label{eq:actionNL}
J^\e(t)\( \(|x|^{-\gamma}\ast|\phi|^2 \)\phi\) =
\(|x|^{-\gamma}\ast|\phi|^2 \)J^\e(t)\phi +2 \operatorname{Re}
\(|x|^{-\gamma}\ast \(\overline{\phi}J^\e(t)\phi\) \)\phi\, .
\end{equation}
The same holds for $H^\e(t)$.
\end{lem}
\begin{rema*}
Property \eqref{eq:commut} follows from \eqref{eq:Heisenberg},
which is the way $J^\e$ and $H^\e$ appear in the linear theory
(see e.g. \cite[p.~108]{Thirring}).
Property \eqref{eq:sobolev} is a
consequence of Gagliardo--Nirenberg inequalities and
\eqref{eq:factor}. Finally, \eqref{eq:actionNL} stems from \eqref{eq:factor}.
\end{rema*}

\section{``Very weak nonlinearity''  case}
\label{sec:small}
In this section, we study the semi-classical limit of $u^\e$ when
{\bf $\gamma =1$ and $\alpha >1$}
which is equivalent to ``very small'' data in our context
(cf. \eqref{eq:r3.1}).
This case includes the 3D Schr\"odinger--Poisson equation with
``very small data''.
We prove that the Hartree term plays no role at
leading order.
\begin{prop}\label{prop:small}
Let $f\in \Sigma$, $n\geq 2$, and assume $\alpha>\gamma=1$. Then for
any $T>0$,
\begin{equation*}
\left\|u^\e-u^\e_{\rm
free}\right\|_{L^\infty([0,T];L^2)}
=\O\(\e^{\alpha-1}\ln\frac{1}{\e}\) \quad \text{as }\e \to 0\, ,
\end{equation*}
and for any $\delta >0$ ($\alpha -1-\delta>0$),
\begin{equation*}
\left\|A^\e(t)\( u^\e-u^\e_{\rm
free}\)\right\|_{L^\infty([0,T];L^2)}
=\O\(\e^{\alpha-1-\delta}\) \quad \text{as }\e \to 0\, ,
\end{equation*}
where $A^\e$ is either of the operators $J^\e$ or
$H^\e$, and $u^\e_{\rm free}$ is the solution of \eqref{eq:harmo}.
\end{prop}
\begin{rema*}
Using modified Sobolev inequalities \eqref{eq:sobolev}, we can deduce
$L^p$ estimates for $u^\e-u^\e_{\rm free}$ for $2\leq p\leq 2n/(n-2)$
($2\leq p<\infty$ if $n=2$) from the above result.
\end{rema*}
\begin{rema*}
We could probably get the logarithmic estimate for the second part of the
statement as well, using Strichartz estimates. The proof given below
is not technically involved, and suffices for our purpose: we do not seek
sharp results.
\end{rema*}
\begin{proof}
Denote $w^\e = u^\e-u^\e_{\rm free}$. It solves the initial value
problem
\begin{equation*}
i\e \d_t w^\e +\frac{1}{2}\e^2\Delta w^\e = \frac{|x|^2}{2}w^\e +
\e^{\alpha}\(|x|^{-1}\ast |u^\e|^2\)u^\e \quad ; \quad w^\e_{\mid t=0} =
0\, .
\end{equation*}
Standard energy estimates for Schr\"odinger equations yield
\begin{equation}\label{14:51}
\e \d_t \|w^\e(t)\|_{L^2} \lesssim \e^\alpha \left\|\(|x|^{-1}\ast
|u^\e|^2\)u^\e \right\|_{L^2}\, .
\end{equation}
From H\"older's inequality, we have
\begin{equation}\label{14:45}
\left\|\(|x|^{-1}\ast
|u^\e|^2\)u^\e \right\|_{L^2}\leq \left\||x|^{-1}\ast
|u^\e|^2\right\|_{L^r}\left\|u^\e \right\|_{L^k}\, ,
\text{ for }\frac{1}{r}
+\frac{1}{k}=\frac{1}{2}\, .
\end{equation}
From the Hardy--Littlewood--Sobolev inequality,
\begin{equation}
  \label{eq:14:47}
 \left\||x|^{-1}\ast |u^\e(t)|^2\right\|_{L^r}
\lesssim \|u^\e(t)\|^2_{L^p}\, ,
\text{ for }1<r,\frac{p}{2}<\infty \text{ and }1+\frac{1}{r}=
\frac{2}{p}+\frac{1}{n}\, .
\end{equation}
Therefore, \eqref{14:51} yields
\begin{equation}\label{eq:15:02}
 \e \d_t \|w^\e(t)\|_{L^2} \lesssim \e^\alpha  \|u^\e(t)\|^2_{L^p}
\left\|u^\e(t) \right\|_{L^k}\, ,
\end{equation}
where $p$ and $k$ satisfy the properties stated in \eqref{14:45} and
\eqref{eq:14:47}. For $k=2$,
$r=\infty$ and $p=2n/(n-1)$, the algebraic identities
 stated in \eqref{14:45} and \eqref{eq:14:47} are satisfied.
Now since the conditions
$1<r<\infty$ and $1<p/2<\infty$ are open, a continuity argument shows that
we can find $p$ and $k$ satisfying
all the properties stated in \eqref{14:45} and
\eqref{eq:14:47}. Notice that they imply the relation
$2\delta(p) +\delta(k)=1$, hence $\delta(p)$, $\delta(k)<1$; this allows us to
use weighted Gagliardo--Nirenberg inequalities.

We have $w^\e_{\mid t=0} =0$, and from Proposition~\ref{prop:cauchy},
$w^\e \in C(\R_+;\Sigma)$. Therefore, there exists $t^\e>0$ such that
\begin{equation}\label{eq:solong.1}
\|J^\e(t)w^\e\|_{L^2}\leq 1\,,
\end{equation}
for $0\leq t\leq t^\e$. The argument of the proof then follows
\cite{RauchUtah} (see also \cite{CaIHP}). Recall that from
Lemma~\ref{lem:operators}, $\|J^\e(t)u^\e_{\rm
free}\|_{L^2}=\|\nabla f\|_{L^2}$.

Because of \eqref{eq:commut},
$J^\e u^\e_{\rm free}$ solves the linear
Schr\"odinger equation with harmonic potential, and
$\|J^\e(t)u^\e_{\rm free}\|_{L^2}\equiv \|\nabla f\|_{L^2}$.
So long as \eqref{eq:solong.1} holds,
we have, from \eqref{eq:sobolev},
\begin{equation*}
 \|u^\e(t)\|_{L^p}\leq \frac{C_0}{|\cos t|^{\delta(p)}}\quad ;\quad
 \|u^\e(t)\|_{L^k}\leq \frac{C_0}{|\cos t|^{\delta(k)}} \, ,
\end{equation*}
for some $C_0$ independent of $\e$ and $t$. Then \eqref{eq:15:02} yields
\begin{equation*}
 \e \d_t \|w^\e(t)\|_{L^2} \lesssim \frac{\e^\alpha }{|\cos t|^{2\delta(p)+
\delta(k)}}= \frac{\e^\alpha}{|\cos t|}\, .
\end{equation*}
Integration in time on $[0,t]$ yields, so long as \eqref{eq:solong.1} holds,
\begin{equation*}
\|w^\e\|_{L^\infty([0,t];L^2)} \lesssim \e^{\alpha
-1} \int_0^t
 \frac{d\tau}{|\cos \tau|}\, ,
\end{equation*}
For $t<\pi/2$, we get, so long as
\eqref{eq:solong.1} holds,
\begin{equation*}
\|w^\e\|_{L^\infty([0,t];L^2)}\lesssim \e^{\alpha -1} \left| \ln
\(\frac{\pi}{2}-t\)\right| \, .
\end{equation*}
From \eqref{eq:commut}, $J^\e(t)w^\e$ solves
\begin{equation*}
i\e \d_t J^\e w^\e +\frac{1}{2}\e^2\Delta J^\e w^\e =
\frac{|x|^2}{2}J^\e w^\e +
\e^{\alpha}J^\e\(\(|x|^{-1}\ast |u^\e|^2\)u^\e\) \quad ; \quad J^\e
w^\e_{\mid t=0} =  0\, .
\end{equation*}
Using \eqref{eq:actionNL}, energy estimate for $J^\e w^\e$ yields
\begin{equation*}
\begin{aligned}
\e \d_t \|J^\e (t)w^\e\|_{L^2}& \lesssim \e^\alpha \(
\left\|\(|x|^{-1}\ast
|u^\e|^2\)J^\e (t)u^\e \right\|_{L^2} + \left\||x|^{-1}\ast
\(\overline{u^\e}J^\e u^\e\) \cdot u^\e \right\|_{L^2}   \) \\
&\lesssim \e^\alpha \(
\left\||x|^{-1}\ast
|u^\e|^2\right\|_{L^\infty}\|J^\e (t) u^\e\|_{L^2} +
\left\||x|^{-1}\ast
\(\overline{u^\e}J^\e u^\e\) \cdot u^\e \right\|_{L^2}   \)
\end{aligned}
\end{equation*}
For the first term of the right hand side, use the easy estimate
\begin{equation*}
  \left\| |x|^{-1}\ast f\right\| \lesssim \|f\|_{L^{(n^-)'}} +
\|f\|_{L^{(n^+)'}}\,
\end{equation*}
where $n^-$ (res. $n^+$) stands for $n -\eta$ (resp. $n+\eta$) for any
small $\eta>0$. We have
\begin{equation*}
 \left\||x|^{-1}\ast
|u^\e|^2\right\|_{L^\infty} \lesssim \| u^\e(t)\|^2_{L^{\kappa^-}} +
\| u^\e(t)\|^2_{L^{\kappa^+}} \, ,
\quad \text{with}\quad \kappa=\frac{2n}{n-1}\, \cdot
\end{equation*}
It is at this stage that we lose the logarithmic rate (we cannot
use Hardy--Littlewood--Sobolev inequality when an exponent is infinite):
using Strichartz estimates (see Section~\ref{sec:scattering}), we believe
that we could recover that rate, with a more technically involved proof.

For the second term, we proceed as in the beginning of the proof. From
H\"older's inequality,
\begin{equation}
  \label{eq:17:50}
  \left\|\(|x|^{-1}\ast
\overline{u^\e}J^\e u^\e\) \cdot u^\e \right\|_{L^2} \leq \left\| |x|^{-1}\ast
\(\overline{u^\e}J^\e u^\e\)\right\|_{L^r} \|u^\e\|_{L^\sigma}\, ,
\text{ with }\frac{1}{r}+\frac{1}{\sigma}=\frac{1}{2}\, \cdot
\end{equation}
From the Hardy--Littlewood--Sobolev inequality, this is estimated, up to a
 constant, by
\begin{equation}
  \label{eq:17:55}
\left\|\overline{u^\e}J^\e u^\e\right\|_{L^p} \|u^\e\|_{L^\sigma}\, ,
\text{ with }1+\frac{1}{r}=\frac{1}{p}+\frac{1}{n}\, \,\text{ for }
1<r,p<\infty\, .
\end{equation}
Using H\"older's inequality again yields an estimate by
\begin{equation}
  \label{eq:17:58}
\left\|u^\e\right\|_{L^k}
\left\|J^\e u^\e\right\|_{L^2} \|u^\e\|_{L^\sigma}\, ,
\text{ with }\frac{1}{p}= \frac{1}{2} + \frac{1}{k} \, \cdot
\end{equation}
Take $r=n$, $\sigma = 2n/(n-2)$, $k=2$ and $p=1$: the algebraic identities
from \eqref{eq:17:50}, \eqref{eq:17:55} and \eqref{eq:17:58} are satisfied,
but not the bound $p>1$. Decreasing slightly $\sigma$ increases $p$
(take $\sigma$ large but finite when $n=2$), so
we can find indices satisfying \eqref{eq:17:50}, \eqref{eq:17:55} and
\eqref{eq:17:58} by a continuity argument. Note that they satisfy
$\delta(k)+\delta(\sigma)=1$, and each term is positive.

Gathering all these estimates together yields the energy estimate
\begin{equation*}
  \e \d_t \|J^\e (t)w^\e\|_{L^2} \lesssim \e^\alpha \(
\| u^\e(t)\|^2_{L^{\kappa^-}} +
\| u^\e(t)\|^2_{L^{\kappa^+}}+ \left\|u^\e\right\|_{L^k}
 \|u^\e\|_{L^\sigma}
\)\|J^\e(t)u^\e\|_{L^2}
\end{equation*}
So long as \eqref{eq:solong.1} holds, we deduce from \eqref{eq:sobolev},
\begin{equation*}
\begin{aligned}
  \e \d_t \|J^\e (t)w^\e\|_{L^2} &\lesssim \e^\alpha \(
\frac{1}{|\cos t|^{2\delta({\kappa^-})}} +
\frac{1}{|\cos t|^{2\delta({\kappa^+})}} +
\frac{1}{|\cos t|^{\delta(k)+\delta(\sigma)}}
\)\\
&\lesssim \e^\alpha \(\frac{1}{|\cos t|^{2\delta({\kappa^+})}} +
\frac{1}{|\cos t|}\)\lesssim  \frac{\e^\alpha}{|\cos t|^{1^+}}\, .
\end{aligned}
\end{equation*}
Integrate this, so long
as \eqref{eq:solong.1} holds:
\begin{equation*}
  \|J^\e w^\e\|_{L^\infty([0,t];L^2)} \lesssim \e^{\alpha -1} \(
\frac{\pi}{2} -t\)^{0^-}
\end{equation*}

Fix $\delta,\Lambda >0$. So long as \eqref{eq:solong.1} holds, we infer, for
$t\leq \pi/2-\Lambda \e$,
\begin{equation*}
\|J^\e w^\e\|_{L^\infty([0,t];L^2)} \lesssim \e^{\alpha-1}
\(\Lambda \e\)^{-\delta} \, .
\end{equation*}
Therefore, there exists $\e_\Lambda >0$ such that, for $0<\e\leq
\e_\Lambda$, \eqref{eq:solong.1} holds up to time $\pi/2-\Lambda \e$,
with the estimates
\begin{equation}\label{eq:est.1}
\|w^\e\|_{L^\infty([0,\pi/2-\Lambda \e];L^2)}\lesssim \e^{\alpha-1} \ln
\frac{1}{\e}\quad ;\quad
\|J^\e
w^\e\|_{L^\infty([0,\pi/2-\Lambda \e];L^2)}\lesssim
\e^{\alpha-1-\delta}  \, .
\end{equation}
An estimate similar to that of $J^\e w^\e$ then follows for $H^\e w^\e$,
since from \eqref{eq:conservations},
$\|H^\e(t)u^\e\|_{L^2}\lesssim \|f\|_\Sigma$. \\

Denote $I^\e_\Lambda = [\pi/2 -\Lambda \e,\pi/2+ \Lambda \e]$. Mimicking the
above computations, we have
\begin{equation*}
\|w^\e\|_{L^\infty(I^\e_\Lambda;L^2)}
\lesssim \left\|w^\e\( \frac{\pi}{2}-\Lambda
\e\)\right\|_{L^2}+ \e^{\alpha
-1}\int_{I^\e_\Lambda}\|u^\e(\tau)\|^2_{L^p}\|u^\e(\tau)\|_{L^k}d\tau\, ,
\end{equation*}
where $p$ and $k$ satisfy \eqref{14:45} and \eqref{eq:14:47}. Recall
that they satisfy $2\delta(p)+\delta(k)=1$.
Using the conservations of mass and energy \eqref{eq:conservations},
along with Gagliardo--Nirenberg inequalities, we have, for any $t$,
\begin{equation*}
\|u^\e(t)\|_{L^p}\lesssim \e^{-\delta(p)}\quad; \quad
\|u^\e(t)\|_{L^k}\lesssim \e^{-\delta(k)}\, .
\end{equation*}
We deduce
\begin{equation*}
\|w^\e\|_{L^\infty(I^\e_\Lambda;L^2)}
\lesssim \left\|w^\e\( \frac{\pi}{2}-\Lambda
\e\)\right\|_{L^2}+ \e^{\alpha
-1}\e^{-2\delta(p)-\delta(k) }\left|I^\e_\Lambda\right| \lesssim
\e^{\alpha-1} \ln
\frac{1}{\e}+ \Lambda \e^{\alpha
-1} \, .
\end{equation*}
The same method yields, since
\eqref{eq:conservations} shows that $\|H^\e (t)u^\e\|_{L^2}\lesssim
\|f\|_\Sigma$:
\begin{equation*}
\|H^\e w^\e\|_{L^\infty(I^\e_\Lambda;L^2)}
\lesssim \e^{\alpha -1-\delta}\, ,\quad \text{for any }\delta>0\, .
\end{equation*}
To treat the case of $J^\e w^\e$, introduce
\begin{equation*}
z_\e(t)= \sup_{\frac{\pi}{2}-\Lambda\e\leq \tau \leq t}
\left\| J^\e(\tau)w^\e\right\|_{L^2}\, .
\end{equation*}
Proceeding as above, we have
\begin{equation}
\label{eq:gron}
\begin{aligned}
z_\e(t) &\lesssim  \left\| J^\e\(\frac{\pi}{2}-\Lambda \e
\)w^\e
\right\|_{L^2}+ \e^{\alpha-1} \int_{\frac{\pi}{2}-\Lambda \e}^t \left\|
J^\e(\tau)\( |x|^{-1}\ast|u^\e|^2 u^\e\)\right\|_{L^2} d\tau \\
&\lesssim \e^{\alpha -1^+} + \e^{\alpha
-1}\int_{\frac{\pi}{2}-\Lambda \e}^t \e^{-1^+} \(
z_\e(\tau) + \left\|
J^\e(\tau)u^\e_{\rm free}\right\|_{L^2} \)d\tau \, .
\end{aligned}
\end{equation}
We can then apply the Gronwall lemma (recall that $\|
J^\e(\tau)u^\e_{\rm free}\|_{L^2}\equiv \|\nabla f\|_{L^2}$):
\begin{equation*}
z_\e(t)
\lesssim \e^{\alpha -1^+}  \, .
\end{equation*}
Gathering these informations we get, for any $\delta>0$:
\begin{align*}
\|w^\e\|_{L^\infty(I^\e_\Lambda;L^2)}& \lesssim  \e^{\alpha-1} \ln
\frac{1}{\e}\, ,\\
 \|J^\e
w^\e\|_{L^\infty(I^\e_\Lambda;L^2)} + \|H^\e
w^\e\|_{L^\infty(I^\e_\Lambda;L^2)}&\lesssim \e^{\alpha -1 -\delta}\, .
\end{align*}
For $t\in [\pi/2+\e ,\pi]$, we can use the same proof as
for $t\in [0,\pi/2-\e]$, to obtain:
\begin{align*}
\|w^\e\|_{L^\infty([0,\pi];L^2)}& \lesssim  \e^{\alpha-1} \ln
\frac{1}{\e}\, ,\\
 \|J^\e
w^\e\|_{L^\infty([0,\pi];L^2)} + \|H^\e
w^\e\|_{L^\infty([0,\pi];L^2)}&\lesssim \e^{\alpha -1 -\delta}\, .
\end{align*}
Repeating the same argument a finite number of times covers any given
time interval $[0,T]$ and completes the proof of
Proposition~\ref{prop:small}.
\end{proof}


\section{Nonlinear propagation and linear focus}
\label{sec:nlwkb}
In this paragraph, we assume $\alpha =1$ and $\gamma <1$. We define
\begin{equation}\label{eq:g}
g(t,x) =- \(|x|^{-\gamma}\ast |f|^2\)(x) \int_0^t \frac{d\tau}{|\cos
\tau|^\gamma}\, .
\end{equation}
This function is well defined for any $t$, since $\gamma <1$. We will
see later on how this function appears.
\begin{prop}\label{prop:nlwkb}
Let $n\geq 2$, $f\in\Sigma$, and assume $\gamma <\alpha =1$. Let
$A^\e$ be either of the operators $Id$, $J^\e$ or
$H^\e$.\\
$\bullet$ For $0\leq t <\pi/2$, the following asymptotics holds:
\begin{equation*}
\sup_{0\leq \tau\leq t} \left\|A^\e(\tau)\( u^\e(\tau ,x)-
\frac{1}{(\cos \tau)^{n/2}}
f\(\frac{x}{\cos\tau}\) e^{-i\frac{x^2}{2\e}\tan \tau +ig\(
\tau,\frac{x}{\cos\tau}\)} \)\right\|_{L^2_x}
\Tend \e 0 0 \, .
\end{equation*}
$\bullet$ For $\pi/2< t \leq\pi$,
\begin{equation*}
\sup_{t\leq \tau\leq \pi} \left\|A^\e(\tau)\( u^\e(\tau ,x)-
\frac{e^{-in\frac{\pi}{2}}}{(\cos \tau)^{n/2}}
f\(\frac{x}{\cos\tau}\) e^{-i\frac{x^2}{2\e}\tan \tau +ig\(
\tau,\frac{x}{\cos\tau}\)} \)\right\|_{L^2_x}
\Tend \e 0 0 \, .
\end{equation*}
$\bullet$ For $t =\pi/2$,
\begin{equation*}
\left\|B^\e\( u^\e\(\frac{\pi}{2}\)-
\frac{1}{\e^{n/2}}\F\(
fe^{ig\(\frac{\pi}{2}\)} \)\(\frac{\cdot}{\e}\)\)\right\|_{L^2}
\Tend \e 0 0 \, ,
\end{equation*}
where $B^\e$ is either of the operators $Id$, $\frac{x}{\e}$ or
$\e\nabla_x$, and the Fourier transform is defined by
\begin{equation}\label{eq:Fourier}
\F \phi (\xi)= \widehat{ \phi (\xi)}=\frac{1}{(2\pi)^{n/2}}\int_{\R^n}e^{-ix\cdot
\xi}\phi(x)dx\, .
\end{equation}
\end{prop}
\begin{rema*}
We can also prove estimates for arbitrarily large time intervals, with
the same proof as below.
\end{rema*}
\begin{rema*}
The difference between the asymptotics before and after the focus is
measured only by the Maslov index, through the phase shift
$e^{-in\pi/2}$: no nonlinear phenomenon occurs at leading order near
the focus. On the other hand, nonlinear effects are relevant outside
the focus, as shown by the presence of $g$.
\end{rema*}

\subsection{Oscillatory integrals}\label{sec:osc}
The main tool for proving Proposition~\ref{prop:nlwkb} is the same as
in linear cases (\cite{Du}, see also \cite{JMRMemoir,Ca2} for
applications in  nonlinear settings): we represent
the solution $u^\e$ as an oscillatory integral. Recall that $u^\e\in
C(\R;\Sigma)$ and that $e^{-i\frac{t}{2\e}(-\e^2\Delta+x^2)}=\U(t)$ is
a unitary group on $L^2$. Define $a^\e$ by
\begin{equation}\label{eq:defa}
a^\e(t,x)=\U(-t)u^\e(t,x)\, .
\end{equation}
We first seek a limit as $\e\to 0$ for $a^\e$ before the focus. This is
suggested by a formal computation as in \cite{CaCMP}, and the
following lemma:
\begin{lem}\label{lem:DA.1}
For $t\in [0,\pi/2[\cup ]\pi/2,\pi]$, define $\V$ by
\begin{equation}\label{eq:linapp}
\V(t)\phi(x) =
\left\{
\begin{aligned}
\frac{1}{(\cos t)^{n/2}}
\phi\(\frac{x}{\cos t}\) e^{-i\frac{x^2}{2\e}\tan t } & &\text{ if
}0\leq t<\pi/2\, ,\\
 \frac{e^{-in\pi/2}}{|\cos t|^{n/2}}
\phi\(\frac{x}{\cos t}\) e^{-i\frac{x^2}{2\e}\tan t } & &\text{ if
}\pi/2< t\leq \pi\, .
\end{aligned}
\right.
\end{equation}
For any $\phi\in
H^1(\R^n)$, any $\theta\in ]0,1/2]$, and any $t\in [0,\pi/2[\cup
]\pi/2,\pi]$,
\begin{equation*}
\left\| \U(t)\phi -\V(t)\phi\right\|_{L^2}\leq 2 |\e \tan t|^\theta
\|\phi\|_{H^1}\, .
\end{equation*}
\end{lem}
\begin{proof}
Notice that from Mehler's formula \eqref{eq:mehler}, we can write, for
$0<t<\pi$,
\begin{equation*}
\U(t)=\M_t \D_t \F \M_t \quad \text{where}\quad
 \M_t(x) = e^{-i\frac{x^2}{2\e\tan t}}\ , \ \D_t
\phi(x)=\frac{1}{(i\e\sin t)^{n/2}}\phi\(\frac{x}{\sin t}\)\, ,
\end{equation*}
and the Fourier transform is defined by \eqref{eq:Fourier}. We infer
\begin{equation*}
\left\| \U(t)\phi -\V(t)\phi\right\|_{L^2}= \left\| \frac{1}{(2i\pi
\tan t)^{n/2}} \int e^{i\frac{|x-y|^2}{2\e \tan t}}f(y)dy -
f(x)\right\|_{L^2}
\end{equation*}
From Parseval formula,
\begin{equation*}
 \frac{1}{(2i\pi
\tan t)^{n/2}} \int e^{i\frac{|x-y|^2}{2\e \tan t}}f(y)dy =
\frac{1}{(2\pi)^{n/2}} \int e^{-i\e \tan t \frac{\xi^2}{2} +ix\cdot
\xi} \F f(\xi)d\xi\, ,
\end{equation*}
therefore
\begin{equation*}
\begin{aligned}
\left\| \U(t)\phi -\V(t)\phi\right\|_{L^2}&= \frac{1}{(2\pi)^{n/2}}
\left\| \int \(e^{-i\e \tan t \frac{\xi^2}{2}}-1\)e^{ ix\cdot
\xi} \F f(\xi)d\xi\right\|_{L^2}\\
&= \left\| \(e^{-i\e \tan t \frac{\xi^2}{2}}-1\) \F
f(\xi)\right\|_{L^2}\, ,
\end{aligned}
\end{equation*}
from Plancherel formula. The lemma then follows from the estimate
$|e^{is}-1|\leq 2 |s|^\theta$, for $0\leq \theta \leq 1/2$.
\end{proof}
From Duhamel's principle, we have
\begin{equation*}
u^\e(t)=\U(t)f -i\int_0^t \U(t-s)\(\(|x|^{-\gamma}\ast
|u^\e|^2\)u^\e\)(s)ds\, .
\end{equation*}
Using \eqref{eq:defa}, we deduce
\begin{equation}\label{eq:dta}
\d_t a^\e(t) = -i \,\U(-t)\(\(|x|^{-\gamma}\ast
|u^\e|^2\)u^\e\)(t)\, .
\end{equation}
Now the formal computation begins. Assume $a^\e\to a$ as $\e\to 0$, in
some suitable sense. Then $u^\e(t)\sim \U(t)a(t)$, and from
Lemma~\ref{lem:DA.1},
\begin{equation*}
u^\e(t,x)\Eq \e 0 \frac{1}{(\cos t)^{n/2}} a\( t,\frac{x}{\cos t}\)
e^{-i\frac{x^2}{2\e}\tan t } \quad \text{for}\quad 0\leq t<\pi/2\, .
\end{equation*}
Plugging this into \eqref{eq:dta} and using Lemma~\ref{lem:DA.1} again
(with $\U(-t)$ instead of $\U(t)$, the result still holds) yields
\begin{equation*}
\d_t a (t,x) = \frac{-i}{|\cos t|^\gamma}\( |x|^{-\gamma}\ast |a|^2\)
a(t,x)\, .
\end{equation*}
Recall that $a_{\mid t=0}= u^\e_{\mid t=0}=f$, and notice that from
the above ordinary differential equation, $\d_t |a|^2=0$: we have
$a(t,x) = f(x)e^{ig(t,x)}$, where
\begin{equation*}
\d_t g (t,x) =\frac{-1}{|\cos t|^\gamma}\(|x|^{-\gamma}\ast
|f|^2\)(x)\quad ; \quad g_{\mid t=0}=0\, .
\end{equation*}
Integrating this equation yields the definition of $g(t,x)$ given in  \eqref{eq:g}.

Proposition~\ref{prop:nlwkb} stems from the more precise following
proposition, Lemma~\ref{lem:DA.1} and a density argument.  In view
of a rigorous justification, denote
\begin{equation}\label{eq:be}
b^\e(t,x) = a^\e(t,x)e^{-ig(t,x)} = e^{-ig(t,x)} \U(-t)u^\e(t,x)\,.
\end{equation}
\begin{prop}\label{prop:nlwkb.1}
Let $f\in \Sigma \cap H^2(\R^n)$. Fix $\delta >0$. There exists
$C_\delta$ such that
\begin{equation*}
\sup_{0\leq t\leq \pi}\|b^\e(t) -f\|_{\Sigma} \leq \int_0^\pi \| \d_t
b^\e(t)\|_{\Sigma} dt \leq C_\delta \e^{1-\gamma -\delta}\, .
\end{equation*}
\end{prop}
The first inequality is trivial. We prove the second one in three
steps:
\begin{itemize}
\item[(i)] On $[0,\pi/2-\Lambda \e]$ for any $\Lambda >0$,
with a constant depending on $\delta$ and $\Lambda$.
\item[(ii)] On $[\pi/2-\Lambda \e,\pi/2+\Lambda \e]$,
with a constant depending on $\delta$ and $\Lambda$.
\item[(iii)] On $[\pi/2+\Lambda \e,\pi]$,
with a constant depending on $\delta$ and $\Lambda$.
\end{itemize}
As in Section~\ref{sec:small}, the parameter $\Lambda>0$ is
arbitrary, while it has to be large in the case $\alpha =\gamma>1$ (see
Section~\ref{sec:scattering} and \cite{CaIHP}): this situation is
typical from a case where the focus is ``linear''.

\subsection{Asymptotics before the focus}\label{sec:avant}
Fix $\Lambda,\delta >0$. We prove that there exists
$C_{\Lambda,\delta}$ such that
\begin{equation}\label{eq:error.1}
\int_0^{\frac{\pi}{2}-\Lambda \e} \| \d_t
b^\e(t)\|_{\Sigma} dt \leq C_{\Lambda,\delta} \e^{1-\gamma
-\delta}\, .
\end{equation}
Denote
\begin{equation*}
y_\e(t)= \int_0^t \| \d_t
b^\e(\tau)\|_{H^1} d\tau\,.
\end{equation*}
From \eqref{eq:dta} and the definition \eqref{eq:be},
\begin{align}
\|\d_t b^\e(t)\|_{L^2}&= \left\| \U(-t)\(\(|x|^{-\gamma}\ast
|u^\e|^2\)u^\e\)(t)-\frac{1}{|\cos t|^\gamma}\(|x|^{-\gamma}\ast
|f|^2\)a^\e(t)\right\|_{L^2}\nonumber\\
&=\left\| \(|x|^{-\gamma}\ast
|u^\e|^2\)u^\e(t)-\frac{1}{|\cos t|^\gamma}\U(t)\(\(|x|^{-\gamma}\ast
|f|^2\)a^\e\)(t)\right\|_{L^2}\, .\label{eq:dtb.1}
\end{align}
Lemma~\ref{lem:DA.1} suggests that we can replace $\U$ with $\V$ in
the last expression, up to a controllable error. Before going further
into details, we prove two lemmas which will be of constant use in the
proof of Proposition~\ref{prop:nlwkb.1}.
\begin{lem}\label{lem:ind.1}
Assume $\gamma<1$, and let $0<\delta <2(1-\gamma)$. There exist $p$
and $q$ with
\begin{equation*}
2\delta(2p')=\gamma +\frac{\delta}{2} \big(<1\big)\ ,\ p<\frac{n}{\gamma}\,
\quad ;\quad
\delta(2q')=\frac{\gamma +1}{2}+\frac{\delta}{4} \big(<1\big)\ ,\
q<\frac{n}{\gamma +1}\, ,
\end{equation*}
and such that there exists $C$ such that for any $\phi \in C^\infty_c(\R^n)$,
\begin{equation*}
\begin{aligned}
&\left\| |x|^{-\gamma}\ast \phi\right\|_{L^\infty} \leq C \(
\|\phi\|_{L^1} + \|\phi\|_{L^{p'}}\)\, ,\\
&\left\| \nabla\(|x|^{-\gamma}\ast \phi\)\right\|_{L^\infty} \leq C \(
\|\phi\|_{L^1} + \|\phi\|_{L^{q'}}\)\, .
\end{aligned}
\end{equation*}
\end{lem}
\begin{proof}
We have $2\delta(2p')=\gamma$ when $p=n/\gamma$, and
$\delta(2q')=\frac{\gamma +1}{2}$ when $q=n/(\gamma+1)$. Therefore
$p<n/\gamma$ and $ q<n/(\gamma+1)$ if $2\delta(2p')=\gamma
+\delta/2$ and $\delta(2q')=\frac{\gamma +1}{2}+\frac{\delta}{4}$.

Let $\chi \in C^\infty_c (\R_+,[0,1])$ with $\chi \equiv 1$ on
$[0,1]$. We have
\begin{equation*}
\begin{aligned}
\| |x|^{-\gamma}\ast \phi \|_{L^\infty} & \leq \left\| \(\chi
|x|^{-\gamma}\)\ast \phi \right\|_{L^\infty} +
\left\| \((1-\chi)|x|^{-\gamma}\)\ast \phi \right\|_{L^\infty}\\
&\leq \left\| \chi
|x|^{-\gamma}\right\|_{L^p} \| \phi \|_{L^{p'}}
 + \|(1-\chi)|x|^{-\gamma}\|_{L^\infty}\|\phi\|_{L^1} \\
&\leq C\( \| \phi \|_{L^{p'}}  + \|\phi\|_{L^1}\)\, ,
\end{aligned}
\end{equation*}
where we have used $x\mapsto |x|^{-\gamma}\in L^p_{\rm loc}(\R^n)$
because $p<n/\gamma$. The other estimate is similar, since $\nabla
|x|^{-\gamma} = \O(|x|^{-\gamma -1})$.
\end{proof}
\begin{lem}\label{lem:borne}
Let $\gamma <1$ and  $f\in \Sigma \cap H^2(\R^n)$. Recall that $g$ is
defined by \eqref{eq:g}. We have:
\begin{equation*}
|x|^{-\gamma}\ast |f|^2 \in W^{2,\infty}\ ;\ g \in L^\infty_{\rm
 loc}(\R;  W^{2,\infty})\ ; \ fe^{ig}, \(|x|^{-\gamma}\ast
 |f|^2\) fe^{ig} \in L^\infty_{\rm
 loc}(\R;  H^2)\, .
\end{equation*}
\end{lem}
\begin{proof}
From Lemma~\ref{lem:ind.1} and Sobolev embeddings,
\begin{equation*}
\begin{aligned}
\left\||x|^{-\gamma}\ast |f|^2\right\|_{L^\infty}&\lesssim \|f\|_{L^2}^2 +
\|f\|_{L^{2p'}}^2 \lesssim \|f\|_{H^1}^2\, ,\\
\left\|\nabla|x|^{-\gamma}\ast |f|^2\right\|_{L^\infty}&\lesssim
\|f\|_{L^2}^2 +
\|f\|_{L^{2q'}}^2 \lesssim \|f\|_{H^1}^2\, ,\\
\left\|\nabla^2|x|^{-\gamma}\ast |f|^2\right\|_{L^\infty}&\lesssim
\left\|\nabla|x|^{-\gamma}\ast \(\nabla |f|^2\)\right\|_{L^\infty}\\
& \lesssim
\| f\|_{L^2}\|\nabla f\|_{L^2} +
\|f\|_{L^{2q'}}\|\nabla f\|_{L^{2q'}} \lesssim \|f\|_{H^2}^2\, .
\end{aligned}
\end{equation*}
Since $t\mapsto |\cos t|^{-\gamma}\in L^1_{\rm loc}(\R)$, we infer
that $g \in L^\infty_{\rm loc}(\R;  W^{2,\infty})$. The last two
properties follow easily.
\end{proof}

We can now replace $\U$ with $\V$ in \eqref{eq:dtb.1}, up to the
following error. From Lemmas~\ref{lem:DA.1}, \ref{lem:ind.1} and
\ref{lem:borne},
\begin{equation*}
\begin{aligned}
\left\| \( \U(t)-\V(t)\) \(\(|x|^{-\gamma}\ast
|f|^2\)a^\e\)(t)\right\|_{L^2} & \lesssim |\e \tan t|^\theta \left\|
|\(|x|^{-\gamma}\ast
|f|^2\)a^\e(t)\right\|_{H^1}\\
&\lesssim |\e \tan t|^\theta \left\||x|^{-\gamma}\ast
|f|^2\right\|_{W^{1,\infty}}\left\|a^\e(t)\right\|_{H^1}\\
&\lesssim |\e \tan t|^\theta \( \left\| a^\e(t)\right\|_{L^2}+
\left\|\nabla_x a^\e(t)\right\|_{L^2}\)\\
&\lesssim |\e \tan t|^\theta \( \|f\|_{L^2} + \left\|\nabla_x
\(b^\e e^{ig}\)\right\|_{L^2}\) \\
&\lesssim |\e \tan t|^\theta \( 1 + \left\|\nabla_x
b^\e(t)\right\|_{L^2}\) \\
&\lesssim |\e \tan t|^\theta \( 1 + \left\|\nabla_x
\(b^\e(t)-f\)\right\|_{L^2}\) \\
&\lesssim |\e \tan t|^\theta \( 1 + \int_0^t \left\|\d_t b^\e(\tau)
\right\|_{H^1}d\tau\)\, ,
\end{aligned}
\end{equation*}
for $0<\theta \leq 1/2$ to be fixed later. Plugging this estimate into
\eqref{eq:dtb.1} yields
\begin{equation}\label{eq:dtb.2}
\begin{aligned}
\|\d_t b^\e(t)\|_{L^2}
\lesssim & \left\| \(|x|^{-\gamma}\ast
|u^\e|^2\)u^\e(t)-\frac{1}{|\cos t|^\gamma}\V(t)\(\(|x|^{-\gamma}\ast
|f|^2\)a^\e\)(t)\right\|_{L^2}\\
& + \frac{|\e \tan t|^\theta }{|\cos t|^\gamma} (1+y_\e(t))\, .
\end{aligned}
\end{equation}
We check that
\begin{equation}\label{eq:commutV}
\frac{1}{|\cos t|^\gamma}\V(t)\(\(|x|^{-\gamma}\ast
|f|^2\)\phi\) = \(|x|^{-\gamma}\ast
|\V(t)f|^2\)\V(t)\phi\, .
\end{equation}
Since we expect $\V(t)a^\e(t)$ to be close to $\U(t)a^\e(t)=u^\e(t)$
as $\e\to 0$, we estimate the difference
\begin{equation*}
\begin{aligned}
\Big\| \(|x|^{-\gamma}\ast |\V(t)f|^2\) & \( \V(t)a^\e(t) -
\U(t)a^\e(t)\)\Big\|_{L^2}\\
 &\lesssim \left\||x|^{-\gamma}\ast
|\V(t)f|^2 \right\|_{L^\infty} \left\| \(\V(t)- \U(t)\)(b^\e
 e^{ig})\right\|_{L^2} \\
&\lesssim \( \| \V(t)f\|_{L^2}^2 + \| \V(t)f\|_{L^{2p'}}^2\) (\e \tan
t)^\theta \left\| b^\e(t)e^{ig(t)}\right\|_{H^1}\\
&\lesssim \(1 + |\cos t|^{-2\delta (2p')}\) |\e \tan t|^\theta \(
\|b^\e(t)-f\|_{H^1} + \|f\|_{H^1}\)\, ,
\end{aligned}
\end{equation*}
%
using the modified Sobolev inequality \eqref{eq:sobolev}.
Since $2\delta (2p')=\frac{n}{p}>\gamma$, we infer from \eqref{eq:dtb.2} that
\begin{equation}\label{eq:dtb.3}
\begin{aligned}
\|\d_t b^\e(t)\|_{L^2}
&\lesssim  \ \frac{|\e \tan t|^\theta }{|\cos t|^{2\delta (2p')}}
(1+y_\e(t))\\
& + \left\| \(|x|^{-\gamma}\ast
\( |u^\e(t)|^2- |\V(t)f|^2\) \)u^\e(t)\right\|_{L^2}\, .
\end{aligned}
\end{equation}
From Lemma~\ref{lem:ind.1}, the last term is estimated, up to a constant,
by
\begin{equation}\label{eq:estim}
\begin{aligned}
\big\| & |u^\e(t)|^2- |\V(t)f|^2\big\|_{L^1} +  \left\|
|u^\e(t)|^2- |\V(t)f|^2\right\|_{L^{p'}}\lesssim \\
&\lesssim  \left\|  u^\e(t)- \V(t)\(f e^{ig(t)}\)\right\|_{L^2}\(
\left\|  u^\e(t)\right\|_{L^2} + \left\|
\V(t)\(f e^{ig(t)}\)\right\|_{L^2} \)\\
&+ \left\|  u^\e(t)- \V(t)\(f e^{ig(t)}\)\right\|_{L^{2p'}}\(
\left\|  u^\e(t)\right\|_{L^{2p'}} + \left\|
\V(t)\(f e^{ig(t)}\)\right\|_{L^{2p'}} \)\, .
\end{aligned}
\end{equation}
For the first term of the right hand side, we have, since $\U$ is
unitary on $L^2$,
\begin{equation*}
\begin{aligned}
\left\| \U(t)\(b^\e e^{ig}\)- \V(t)\(f e^{ig}\)\right\|_{L^2}
&\lesssim
\|b^\e(t) -f\|_{L^2} + \left\| \( \U(t)-\V(t)\)\(f
e^{ig(t)}\)\right\|_{L^2} \\
& \lesssim y_\e(t) + |\e \tan t|^\theta \left\| f
e^{ig(t)}\right\|_{H^1}\, .
\end{aligned}
\end{equation*}
In addition, notice that $\|u^\e(t)\|_{L^2} = \|\V(t)f\|_{L^2} =
\|f\|_{L^2}$.
The second term is estimated thanks to the modified
Gagliardo--Nirenberg inequality \eqref{eq:sobolev},
\begin{equation*}
\begin{aligned}
\left\|  u^\e(t)- \V(t)\(f e^{ig(t)}\)\right\|_{L^{2p'}} \lesssim & \
|\cos t|^{-\delta(2p')} \left\|  u^\e(t)- \V(t)\(f
 e^{ig(t)}\)\right\|_{L^{2}}^{1-\delta(2p')}\times \\
&\times\left\|J^\e(t)\(
 u^\e(t)- \V(t)\(f
 e^{ig(t)}\)\)\right\|_{L^{2}}^{\delta(2p')} .
\end{aligned}
\end{equation*}
The first $L^2$--norm was estimated just above. For the second one,
notice that
\begin{equation*}
J^\e(t)\U(t) = -i\, \U(t)\nabla_x \quad ; \quad J^\e(t)\V(t) = -i\,
\V(t)\nabla_x  \, ,
\end{equation*}
therefore:
\begin{equation*}
\begin{aligned}
\Big\|J^\e(t) \Big(
 u^\e(t)- &\V(t)\Big(f
 e^{ig(t)}\Big)\Big)\Big\|_{L^{2}}
\lesssim \Big\|\U(t) \nabla\(  b^\e(t) e^{ig(t)}\) - \V(t)\nabla \(f
 e^{ig(t)}\)\Big\|_{L^{2}} \\
&\lesssim \left\|\nabla\(  b^\e(t) e^{ig(t)}- f
 e^{ig(t)}\)\right\|_{L^{2}} + \left\|\(\U(t)-\V(t)\) \nabla\( f
 e^{ig(t)}\)\right\|_{L^{2}}\\
&\lesssim y_\e(t) + |\e \tan t|^\theta\, ,
\end{aligned}
\end{equation*}
where we have used Lemmas~\ref{lem:DA.1} and \ref{lem:borne}. We infer
that
\begin{equation*}
\left\|  u^\e(t)- \V(t)\(f e^{ig(t)}\)\right\|_{L^{2p'}}
\lesssim
|\cos t|^{-\delta(2p')} \( y_\e(t) + |\e \tan t|^\theta\)\, .
\end{equation*}
We have explicitly
\begin{equation*}
\left\| \V(t)\(f e^{ig(t)}\)\right\|_{L^{2p'}} = |\cos
t|^{-\delta(2p')} \|f\|_{L^{2p'}}
\lesssim
|\cos t|^{-\delta(2p')} \, .
\end{equation*}
Proceeding as above, we have
\begin{equation*}
\left\| u^\e(t)\right\|_{L^{2p'}}  \lesssim  |\cos t|^{-\delta(2p')}
\|u^\e\|_{L^2}^{1-\delta(2p')} \|J^\e(t)u^\e\|_{L^2}^{\delta(2p')}\,
,
\end{equation*}
with $ \|J^\e(t)u^\e\|_{L^2} \lesssim \|b^\e(t)-f\|_{H^1} +
\|f\|_{H^1} $. These estimates  will eventually lead to an
inequality of the form $y'_\e(t)\leq a(t)y_\e(t) + b(t)y_\e(t)^\kappa + c(t)$,
for some $\kappa >1$. To avoid that situation, we proceed as in
Section~\ref{sec:small}; there exists $t^\e>0$ such that
\begin{equation}\label{eq:solong.2}
\|b^\e(t)\|_{H^1}\leq 2\|f\|_{H^1}\, ,
\end{equation}
for $t\in [0,t^\e]$. So long as \eqref{eq:solong.2} holds, we have
from the above estimates
\begin{equation}\label{eq:dtb.4}
\|\d_t b^\e(t)\|_{L^2} \lesssim |\cos t|^{-2\delta(2p')}\( y_\e(t) + |\e
\tan t|^\theta\)\, .
\end{equation}
To prove that \eqref{eq:solong.2} holds up to time $\pi/2 -\Lambda \e$
for $0<\e\leq \e_\Lambda$ along with the error estimate
\eqref{eq:error.1}, we estimate the $L^2$--norm of $\nabla_x \d_t
b^\e$.  From \eqref{eq:dta} and \eqref{eq:be},
\begin{align*}
\nabla_x \d_t b^\e(t)&= -i\nabla_x g(t) \d_t b^\e(t)\\
-i e^{-ig(t)} & \nabla_x \(\U(-t)\(\(|x|^{-\gamma}\ast
|u^\e|^2\)u^\e\)(t)-\frac{1}{|\cos t|^\gamma}\(|x|^{-\gamma}\ast
|f|^2\)a^\e(t)\)\, .
\end{align*}
The first term is controlled thanks to Lemma~\ref{lem:borne} and
\eqref{eq:dtb.4}. For the other term, we notice that since $\U$ is
unitary on $L^2$, from \eqref{eq:Heisenberg} its $L^2$--norm is equal
to:
\begin{equation*}
\left\| J^\e(t)\(\(|x|^{-\gamma}\ast
|u^\e|^2\)u^\e\)(t)+\frac{i}{|\cos
t|^\gamma}\U(t)\nabla_x\(\(|x|^{-\gamma}\ast
|f|^2\)a^\e(t)\)\right\|_{L^2}\, .
\end{equation*}
We proceed as before: we first replace $\U$ with $\V$ in the last
term, up to an error of $|\cos t|^{-\gamma}$ times:
\begin{align*}
\Big\| \(\U(t)-\V(t)\)&\nabla_x\(\(|x|^{-\gamma}\ast
|f|^2\)a^\e\)\Big\|_{L^2}\lesssim \\
&\lesssim\left\|
\(\U(t)-\V(t)\)\nabla_x\(\(|x|^{-\gamma}\ast
|f|^2\)(b^\e -f)e^{ig}\)\right\|_{L^2}\\
&+\left\|
\(\U(t)-\V(t)\)\nabla_x\(\(|x|^{-\gamma}\ast
|f|^2\) fe^{ig}\)\right\|_{L^2} \, .
\end{align*}
For the first term, we do not use Lemma~\ref{lem:DA.1}, but roughly
the fact that $\U$ and $\V$ are unitary on $L^2$. It is not larger than
\begin{equation*}
2\left\|\nabla_x\(\(|x|^{-\gamma}\ast
|f|^2\)(b^\e -f)e^{ig}\)\right\|_{L^2}\lesssim \|b^\e(t)-f\|_{H^1}\, ,
\end{equation*}
from Lemma~\ref{lem:borne}. The second term is controlled thanks to
Lemmas~\ref{lem:DA.1} and \ref{lem:borne},
\begin{equation*}
\left\|
\(\U(t)-\V(t)\)\nabla_x\(\(|x|^{-\gamma}\ast
|f|^2\) fe^{ig}\)\right\|_{L^2} \lesssim  |\e \tan t|^\theta\, .
\end{equation*}
We now have, so long as \eqref{eq:solong.2} holds,
\begin{equation}\label{eq:dtb.5}
\begin{aligned}
&\|\d_t b^\e(t)\|_{H^1} \lesssim |\cos t|^{-2\delta(2p')}\( y_\e(t) + |\e
\tan t|^\theta\) \\
 + \Big\| J^\e(t)&\(\(|x|^{-\gamma}\ast
|u^\e|^2\)u^\e\)+\frac{i}{|\cos
t|^\gamma}\V(t)\nabla_x\(\(|x|^{-\gamma}\ast
|f|^2\)a^\e(t)\)\Big\|_{L^2}\, .
\end{aligned}
\end{equation}
Using the identity $J^\e(t)\V(t) = -i\, \V(t)\nabla_x$ and
\eqref{eq:commutV}, we have to estimate
\begin{equation}\label{eq:inter.1}
\begin{aligned}
&\left\| J^\e(t)\(\(\(|x|^{-\gamma}\ast
|u^\e|^2\)u^\e\)-\frac{1}{|\cos
t|^\gamma}\V(t)\(\(|x|^{-\gamma}\ast
|f|^2\)a^\e(t)\)\)\right\|_{L^2}\\
&=\left\| J^\e(t)\(\(|x|^{-\gamma}\ast
|u^\e|^2\)u^\e-\(|x|^{-\gamma}\ast
|\V(t)f|^2\)\V(t)a^\e\)\right\|_{L^2}\\
&\lesssim  \left\| \(|x|^{-\gamma}\ast
|u^\e|^2\)J^\e(t)u^\e-\(|x|^{-\gamma}\ast
|\V(t)f|^2\)J^\e(t)\V(t)a^\e\right\|_{L^2}\\
&+ |\cos t| \left\| \nabla_x\(|x|^{-\gamma}\ast
|u^\e|^2\) u^\e-\nabla_x\(|x|^{-\gamma}\ast
|\V(t)f|^2\)\V(t)a^\e\right\|_{L^2}\, .
\end{aligned}
\end{equation}
We replace $\V$ with $\U$ in the first term of the right hand side, up
to the error
\begin{align*}
\big\| \(|x|^{-\gamma}\ast
|\V(t)f|^2\)&J^\e(t)\(\V(t)-\U(t)\)a^\e\big\|_{L^2}\lesssim\\
\lesssim  &
\( \|\V(t)f\|_{L^2}^2 +\|\V(t)f\|_{L^{2p'}}^2\)  \left\|
\(\V(t)-\U(t)\)\nabla_x a^\e\right\|_{L^2}\\
\lesssim &|\cos t|^{-2\delta(2p')}  \left\|
\(\V(t)-\U(t)\)\nabla_x \((b^\e-f)e^{ig}\)\right\|_{L^2}\\
&+|\cos t|^{-2\delta(2p')}  \left\|
\(\V(t)-\U(t)\)\nabla_x \((f e^{ig}\)\right\|_{L^2}\\
\lesssim  & |\cos t|^{-2\delta(2p')}  \( \left\|
b^\e(t)-f\right\|_{H^1} + |\e \tan t|^\theta\)\, ,
\end{align*}
from the above computation. Therefore, the first term of the right
hand side of \eqref{eq:inter.1} is estimated by
\begin{equation*}
|\cos t|^{-2\delta(2p')}\(y_\e(t)+|\e \tan t|^\theta\) +
\left\| \(|x|^{-\gamma}\ast
\(|u^\e|^2-|\V(t)f|^2\)\)J^\e(t)u^\e \right\|_{L^2}\, .
\end{equation*}
So long as \eqref{eq:solong.2} holds, $\|J^\e(t)u^\e \|_{L^2}\lesssim
1$, and the last term is estimated by
\begin{equation*}
\left\| |x|^{-\gamma}\ast
\(|u^\e|^2-|\V(t)f|^2\) \right\|_{L^\infty}\, ,
\end{equation*}
which already appeared above and was estimated in \eqref{eq:estim}.
We are left with the second term of the right hand side of
\eqref{eq:inter.1}. Using Lemma~\ref{lem:ind.1} with $q$ instead of
$p$ now,
\begin{align*}
\big\| \nabla_x\(|x|^{-\gamma}\ast
|\V(t)f|^2\)&\(\V(t)-\U(t)\)a^\e\big\|_{L^2}\lesssim \\
&\lesssim
\( \|\V(t)f\|_{L^2}^2 +\|\V(t)f\|_{L^{2q'}}^2  \) \left\|
\(\V(t)-\U(t)\)a^\e\right\|_{L^2}\\
&\lesssim |\cos t|^{-2\delta (2q')} \( y_\e(t) + |\e \tan t |^\theta\)\, .
\end{align*}
The final term to estimate is
\begin{equation*}
\begin{aligned}
\left\| \nabla_x \(|x|^{-\gamma}\ast
\(|u^\e|^2-|\V(t)f|^2\)\)u^\e \right\|_{L^2}\lesssim &
\left\||u^\e|^2-|\V(t)f|^2 \right\|_{L^1}\\
& + \left\||u^\e|^2-|\V(t)f|^2
\right\|_{L^{q'}}\, .
\end{aligned}
\end{equation*}
The right hand side was already estimated in \eqref{eq:estim}
 with $p$ instead of $q$.
We finally have, so long as
\eqref{eq:solong.2} holds,
\begin{equation*}
y'(t) \lesssim \( |\cos t|^{-2\delta(2p')} + |\cos t|^{1-
2\delta(2p')}\)\( y_\e(t) +|\e \tan t|^\theta\)\, .
\end{equation*}
Now recall that given $\delta>0$, $\delta(2p')$ and $\delta(2q')$ are
explicit, hence
\begin{equation*}
y'(t) \lesssim |\cos t|^{-\gamma -\frac{\delta}{2}} \( y_\e(t) +|\e \tan
t|^\theta\)\, .
\end{equation*}
It is now time to fix $\theta$. In view of \eqref{eq:error.1}, it is
natural to take $\theta = 1-\gamma -\delta$. This yields, so long as
\eqref{eq:solong.2} holds,
\begin{equation}\label{eq:gronwall.1}
y'_\e(t) \lesssim |\cos t|^{-\gamma -\frac{\delta}{2}} \( y_\e(t) +|\e \tan
t|^{1-\gamma -\delta}\)\lesssim |\cos t|^{-\gamma -\frac{\delta}{2}}
y_\e(t) + \frac{\e^{1-\gamma -\delta}}{|\cos t|^{1-\frac{\delta}{2}}} \, .
\end{equation}
The maps $t \mapsto |\cos t|^{-\gamma -\frac{\delta}{2}}$ and
$t\mapsto |\cos t|^{-1+\frac{\delta}{2}}$ are locally integrable (we
can assume $\gamma +\delta/2 < 1 -\delta/2$, otherwise
\eqref{eq:error.1} is of no interest). From the Gronwall lemma,  so long as
\eqref{eq:solong.2} holds, we infer
\begin{equation}\label{eq:dtb.6}
y_\e(t) \lesssim \e^{1-\gamma -\delta}\, .
\end{equation}
Therefore, there exists $\e_\Lambda>0$ such that for $0<\e\leq
\e_\Lambda$, \eqref{eq:solong.2} holds up to time $\pi/2 -\Lambda \e$,
with \eqref{eq:dtb.6}. The estimate for $x\d_t b^\e$ then is
easy, we leave out this part; this proves \eqref{eq:error.1}.
\begin{rema*}
One might believe that we could deduce Proposition~\ref{prop:nlwkb.1}
in one shot from \eqref{eq:gronwall.1}, and wonder why we split the
proof into three steps. The reason is that we cannot apply
Lemma~\ref{lem:DA.1} (which was used to get
\eqref{eq:gronwall.1}) near $t=\pi/2$. On the other hand, we will
see below that computations near $t=\pi/2$ are far simpler.
\end{rema*}

\subsection{Near the focus and beyond}
Keep $\Lambda,\delta>0$ fixed. We prove that there exists
$C_{\Lambda,\delta}$ such that
\begin{equation}\label{eq:error.2}
\int_{\frac{\pi}{2}-\Lambda \e}^{\frac{\pi}{2}+\Lambda \e} \| \d_t
b^\e(t)\|_{\Sigma} dt  \leq C_{\Lambda,\delta} \e^{1-\gamma
-\delta}\, .
\end{equation}
A rough estimate in \eqref{eq:dtb.1} yields
\begin{equation}\label{eq:dtb.7}
\begin{aligned}
\|\d_t b^\e(t)\|_{L^2}&\lesssim  \left\| \(|x|^{-\gamma}\ast
|u^\e|^2\)u^\e(t)\right\|_{L^2} +\frac{1}{|\cos t|^\gamma}
\left\|\(|x|^{-\gamma}\ast
|f|^2\)a^\e(t)\right\|_{L^2}\\
&\lesssim \( \|u^\e(t)\|_{L^2}^2 + \|u^\e(t)\|_{L^{2p'}}^2 \)
\|u^\e(t)\|_{L^2} +\frac{1}{|\cos t|^\gamma}
\left\|u^\e(t)\right\|_{L^2}\, .
\end{aligned}
\end{equation}
The conservation of mass yields $\|u^\e(t)\|_{L^2}=\|f\|_{L^2}$. The
conservations of mass and energy \eqref{eq:conservations} yield, along
with Gagliardo--Nirenberg inequalities,
\begin{equation*}
\|u^\e(t)\|_{L^{2p'}} \lesssim \e^{-\delta (2p')}\, .
\end{equation*}
Using this estimate (which is sharp near the focus, and only near the
focus) and integrating \eqref{eq:dtb.7}, we get
\begin{equation*}
\begin{aligned}
\int_{\frac{\pi}{2}-\Lambda \e}^{\frac{\pi}{2}+\Lambda \e} \|\d_t
b^\e(t)\|_{L^2} dt &\lesssim  \Lambda \e^{1-2\delta(2p')} +
\int_{\frac{\pi}{2}-\Lambda \e}^{\frac{\pi}{2}+\Lambda \e}
\frac{dt}{|\cos t|^\gamma} \\
&\lesssim  \e^{1-\gamma-\frac{\delta}{2}} + \e^{1-\gamma}\, .
\end{aligned}
\end{equation*}
The term $\|x\d_t
b^\e(t)\|_{L^2}$ is estimated the same way, since the conservation of
energy yields an \emph{a priori} bound for $H^\e u^\e$. For $\|\nabla_x
\d_t b^\e(t)\|_{L^2}$, we proceed as in Section~\ref{sec:small},
\eqref{eq:gron} to get an estimate from Gronwall lemma; the details are left
to the reader. \\

Finally, one can prove that there exists
$C_{\Lambda,\delta}$ such that
\begin{equation*}
\int_{\frac{\pi}{2}+\Lambda \e}^\pi \| \d_t
b^\e(t)\|_{\Sigma} dt  \leq C_{\Lambda,\delta} \e^{1-\gamma
-\delta}
\end{equation*}
by mimicking the computations performed in Section~\ref{sec:avant},
and the proof of Proposition~\ref{prop:nlwkb.1} is complete.


\section{Linear propagation and nonlinear focus}
\label{sec:scattering}
We now consider the case where $\alpha =\gamma >1$ in
\eqref{eq:r3.3}. Our results are similar to those of
\cite{CaIHP}. Before stating the main result, we recall some points of
the scattering theory for \eqref{eq:r3gen}.
\begin{prop}[\cite{GV80,HT87a}]\label{prop:scattr3}
Assume $\psi_- \in \Sigma$ and $1<\gamma <\min(4,n)$.
If $\gamma>4/3$ or if $\|\psi_-\|_\Sigma$ is
sufficiently small, then:
\begin{itemize}
\item There exists a unique $\psi \in C(\R_t,\Sigma)$ solution of
\eqref{eq:r3gen}, such that
\begin{equation*}
\lim_{t\rightarrow -\infty}\|\psi_- -{\tt U}(-t)\psi(t)\|_\Sigma
=0\, ,\quad\text{where } {\tt U}(t)=e^{i\frac{t}{2}\Delta}\, .
\end{equation*}
\item There exists a unique $\psi_+ \in \Sigma$ such that
\begin{equation*}
\lim_{t\rightarrow +\infty}\|\psi_+ -{\tt U}(-t)\psi(t)\|_\Sigma =0\,
.
\end{equation*}
\end{itemize}
The scattering operator is $S:\psi_-\mapsto \psi_+$.
\end{prop}
Our main result in this section is:
\begin{prop}\label{prop:scattering}
Suppose $n\geq 2$. Let $f\in \Sigma$, $1<\gamma =\alpha<
\min(4,n)$, and $k\in
\N$. Assume either $\gamma >4/3$ or $\|f\|_\Sigma $ is sufficiently
small. Then
the asymptotics of $u^\eps$ for $\pi/2 +
(k-1)\pi < a\leq b < \pi/2 +k\pi$ is given by
\begin{equation*}
\sup_{a\leq t\leq b}\left\| A^\e(t)\(u^\e(t,x) -
\frac{e^{-ink\frac{\pi}{2}}}{|\cos
t|^{n/2}}\(\F\circ S^k\circ\F^{-1}\) f\left(\frac{x}{\cos t} \right)
e^{-i\frac{x^2}{2\e}\tan t}\)\right\|_{L^2_x}\Tend \e 0 0\, ,
\end{equation*}
where $A^\e$ is either of the operators $Id$, $J^\e$ or
$H^\e$, and $S^k$ denotes the $k$-th iterate of $S$ (which is well defined
under our assumptions on $f$). At the focuses:
\begin{equation*}
\left\|B^\e\(u^\e\(\frac{\pi}{2}+k\pi\) -
\frac{e^{-ink\frac{\pi}{2}}}{
\e^{n/2}}\(\F\circ S^k\) f\left(\frac{\cdot}{\e} \right)
\)\right\|_{L^2}\Tend \e 0 0\, ,
\end{equation*}
where $B^\e$ is either of the operators $Id$, $\frac{x}{\e}$ or
$\e\nabla_x$.
\end{prop}
With Lemma~\ref{lem:DA.1} in mind, this shows that nonlinear effects
are negligible away from focuses, while they have an influence at
leading order near the focuses: each caustic crossing is described in
average by the nonlinear scattering operator $S$ (the phase shift
$e^{-ink\frac{\pi}{2}}$ is the Maslov index, present in the linear
case \cite{Du}).

\smallskip

The proof of Proposition~\ref{prop:scattering}
is very similar to the one in \cite{CaIHP}, which relies on (scaled)
Strichartz estimates. We will refrain from repeating everything in detail
and limit ourselves to prove the main technical proposition
and present an outline for the rest of the proof.
One main difference to the problem in  \cite{CaIHP}
is the action of the operators
$J^\e(t)$, $H^\e(t)$ on the Hartree nonlinearity as
described by \eqref{eq:actionNL}.


We start by reformulating Equation~\eqref{eq:r3.3} by the Duhamel formula
\begin{equation}\label{eq:Duhgen}
\begin{split}
u^\eps(t)=\U(t-t_0)u_0^\eps &-i \eps^{\gamma -1}
\int_{t_0}^t \U(t-s)F^\eps (u^\eps)(s) ds \\
&-i\eps^{-1}
 \int_{t_0}^t \U(t-s)h^\eps (s) ds.
\end{split}
\end{equation}
This equation generalizes  Eq.~(\ref{eq:r3.3}) to the case of an
additional source term  and a general nonlinear term
$F^\eps$. The main technical result which is used throughout the proof
of  Proposition~\ref{prop:scattering} is:

\begin{prop}\label{prop:estgen}Let $t_1>t_0$, with $|t_1-t_0|\leq \pi$.
Let ${ q}, { r}, { s},
 { k} \in [1,\infty]$ be such that:
\begin{equation}\label{eq:holder}
  \left\{
  \begin{aligned}
    \mbox{(a)}& \quad  \frac{1}{{r}'}=\frac{1}{{r}}+
\frac{2}{{s}} + \frac{\gamma}{n} - 1
  \quad \mbox{and} \quad  { s}<\frac{2n}{n-\gamma},\\
     \mbox{(b)}& \quad   \frac{1}{{q}'}=\frac{1}{{q}}+
\frac{2}{{k}}, \\
     \mbox{(c)}& \quad
 ( { q}, { r})  \mbox{ is an admissible pair }, \\
     \mbox{(d)}&  \quad
      0 < \frac{1}{{k}} <\delta({s})<1\, . \\
  \end{aligned}\right.
\end{equation}
Assume that there exists a constant $C$ independent of $t$ and $\eps$
such that for $t_0\leq t\leq t_1$,
\begin{equation}\label{eq:hypF}
\|F^\eps(u^\eps)(t)\|_{L^{{ r}'}_x}\leq \frac{C}{\left(|\cos
t|+\eps \right)^{2
\delta({s})}}\|u^\eps(t)\|_{L^{ r}_x}\ \ \ ,
\end{equation}
and define
$$A^\eps(t_0,t_1):= \left(\int_{t_0}^{t_1}
\frac{dt}{\left(|\cos t|+\eps \right)^{{k}
\delta({s})}}\right)^{2/{k}}.$$
Then there exists $C^*$ independent of $\eps$, $t_0$ and $t_1$ such
that for any admissible pair $(\rho,\sigma)$,
\begin{equation}\label{eq:estgen}
\begin{split}
\|u^\eps\|_{L^{ q}(t_0,t_1;L^{ r})} \leq & C^*
\eps^{-1/{ q}}\|u_0^\eps\|_{L^2} + C_{{ q}, \rho}
\eps^{-1-\frac{1}{{ q}}- \frac{1}{\rho}}\|h^\eps\|_
{ L^{\rho'}(t_0,t_1;L^{\sigma'})}\\
& + C^*\eps^{2\left(\delta( s)-\frac{1}{
k}\right)} A^\eps(t_0,t_1) \|u^\eps\|_{L^{
q}(t_0,t_1;L^{ r})},
\end{split}
\end{equation}

\end{prop}

Mostly the following corollary is applied:
\begin{cor}\label{cor:estgen}
Suppose the assumptions of Prop.~\ref{prop:estgen} are
satisfied. Assume moreover that
$C^*\eps^{2\left(\delta( s)-\frac{1}{
k}\right)} A^\eps(t_0,t_1)\leq 1/2$, which holds
in either of the two cases,
\begin{itemize}
\item $0\leq t_0\leq t_1\leq \frac{\pi}{2}-\Lambda \eps$, with
$\Lambda\geq \Lambda_0$ sufficiently large,
\item $t_0,t_1\in[\frac{\pi}{2}-\Lambda \eps,\frac{\pi}{2}+\Lambda
\eps]$, with $\frac{t_1-t_0}{\eps}\leq \eta$ sufficiently small.
\end{itemize}
Then
\begin{equation}\label{eq:estgenL2}
\|u^\eps\|_{L^\infty(t_0,t_1;L^2)} \leq  C
\|u_0^\eps\|_{L^2} + C_{{ q}, \rho}
\eps^{-1- \frac{1}{\rho}}\|h^\eps\|_
{ L^{\rho'}(t_0,t_1;L^{\sigma'})}.
\end{equation}
\end{cor}

\smallskip
To prove Proposition~\ref{prop:estgen}, we first prove the following
algebraic lemma:
\begin{lem}\label{lem:alg}Let $n\geq 2$, and
assume $1<\gamma<\min(4,n)$. Then there exist
${ q}, { r}, { s}, { k} \in [1,\infty]$
satisfying the conditions \eqref{eq:holder}.
\end{lem}
\begin{proof}
Note that (a) is equivalent to demanding $\gamma / 2 = \delta(r)+\delta(s)$
and $\gamma / 2 > \delta(s)$.
\\
{\em Case $\gamma \leq 2$}: Suppose $\gamma / 2 = \delta(s)$.
Then, by the first half of (a)  $\delta(r)=0$ and
$(q,r)=(\infty,2)$ by (c). With $k=2$, (b) and (d) are satisfied.
Now choose $s$ such that
$1/2 < \delta(s) < \gamma / 2$, but close enough to $\gamma/2$ for
 \eqref{eq:holder}
 still to be valid by continuity
( for example $\delta(s) = \frac 1 2
+ \frac 1 2 \(\frac \gamma 2 - \frac 1 2\) $ ). Then \eqref{eq:holder}
is satisfied.\\
{\em Case $\gamma > 2$}: In this case take  $s$ such that $\delta(s) = 1$,
e.g. $s=\frac{2n}{n-2}$. Up to a continuity argument as in the previous
case, $\delta(s)<1$ and \eqref{eq:holder} is satisfied.
\end{proof}

\begin{proof}[Proof of Proposition \ref{prop:estgen}]
Application of the (scaled) Strichartz estimates (Prop.~\ref{prop:strichartz})
to equation \eqref{eq:Duhgen} yields
\begin{equation*}
\begin{split}
\|u^\eps\|_{L^{ q}(t_0,t_1;L^{ r})} \leq & C
\eps^{-1/{ q}}\|u_0^\eps\|_{L^2} + C_{{ q}, \rho}
\eps^{-1-\frac{1}{{ q}}- \frac{1}{\rho}}\|h^\eps\|_
{ L^{\rho'}(t_0,t_1;L^{\sigma'})}\\
& + C\eps^{\gamma -1 -\frac{2}{
q}}\|F^\eps(u^\eps)\|_{L^{{ q}'}(t_0,t_1;L^{{ r}'})}.
\end{split}
\end{equation*}
Then by the assumptions on $F^\eps(u^\eps)$, after an application of
H\"older inequality in time, the statement follows.
\end{proof}

\begin{proof}[Proof of Corollary~\ref{cor:estgen}]
The additional assumption implies that the last term in
\eqref{eq:estgen} can be absorbed by the left hand side, and we get
\begin{equation}\label{eq:abs1}
\|u^\eps\|_{L^{ q}(t_0,t_1;L^{ r})} \leq  C
\eps^{-1/{ q}}\|u_0^\eps\|_{L^2} + C
\eps^{-1-\frac{1}{{ q}}- \frac{1}{\rho}}\|h^\eps\|_
{ L^{\rho'}(t_0,t_1;L^{\sigma'})}.
\end{equation}
Another application of Strichartz estimates to  equation \eqref{eq:Duhgen},
with indices $(\infty,2)$ on the left and
$(\rho,\sigma)$ respectively $( q, r)$ on the right, yields
\begin{equation*}
\begin{split}
\|u^\eps\|_{L^\infty(t_0,t_1;L^2)} \leq & C
\|u_0^\eps\|_{L^2} + C
\eps^{-1- \frac{1}{\rho}}\|h^\eps\|_
{ L^{\rho'}(t_0,t_1;L^{\sigma'})}\\
& + C\eps^{\gamma -1 -\frac{1}{
q}}\|F^\eps(u^\eps)\|_{L^{{ q}'}(t_0,t_1;L^{{ r}'})}.
\end{split}
\end{equation*}
As before,
\begin{equation*}
\begin{split}
\eps^{\gamma -1 -\frac{1}{
q}}\|F^\eps(u^\eps)\|_{L^{{ q}'}(t_0,t_1;L^{{
r}'})}&  \leq C \eps^{\frac{1}{
q}}\eps^{2\left(\delta( s)-\frac{1}{
k}\right)} A^\eps(t_0,t_1) \|u^\eps\|_{L^{
q}(t_0,t_1;L^{ r})}\\
& \leq C \eps^{\frac{1}{
q}}\|u^\eps\|_{L^{
q}(t_0,t_1;L^{ r})},
\end{split}
\end{equation*}
and the statement now follows from \eqref{eq:abs1}.
\end{proof}

\bigskip
The proof of Proposition~\ref{prop:scattering} consists of three parts:
the propagation before the focus, the matching between the two regimes,
and proof that  near the focus, the harmonic potential
is negligible.
In all parts the main tool to derive the major statements
is Prop. \ref{prop:estgen}. Since the proof is very similar to the
one in \cite{CaIHP}, we do not repeat everything in detail
but give a detailed proof only for the first part to show how
the methods of \cite{CaIHP} are applied.

We now show the proof for the propagation before the focus, that is
the approximation of $  u^\eps(t)$ by $u^\e_{\rm free}(t)$ for
$0\leq t \leq \frac{\pi}{2}-\Lambda \eps$, in the limit $\Lambda\to +\infty$. We prove:
\begin{equation*}
\limsup_{\eps \to 0} \sup_{0\leq t\leq \frac \pi 2 - \Lambda \e }
\Big\| A^\e(t)\(u^\e(t,x) - u^\e_{\rm free}(t,x) \)
\Big\|_{L^2_x}\Tend{\Lambda }{+\infty}  0\, ,
\end{equation*}
with $A^\e(t)$ being either of the operators  $Id$, $J^\e$ or $H^\e$.

 Define the remainder $w^\eps = u^\eps - u^\e_{\rm free}$. It solves
\begin{equation*}
\left\{
\begin{split}
i\eps\d_t w^\eps +\frac{1}{2}\eps^2\Delta w^\eps &
=V(x)w^\eps +\eps^\gamma \(|x|^{-\gamma}\ast |u^\e|^2\) u^\e\, ,\\
w^\eps_{\mid t=0}& = 0 \, . 
\end{split}\right.
\end{equation*}
From Duhamel's principle, this can be written as
\begin{equation}\label{eq:DuhW}
w^\eps(t)= \U(t)r^\eps -i \eps^{\gamma-1}\int_0^t
\U(t-s) \(|x|^{-\gamma}\ast |u^\e|^2\) u^\e(s)ds.
\end{equation}
Since $u^\e_{\rm free}$ solves the linear equation (\ref{eq:harmo}), so does
$J^\eps(t)u^\e_{\rm free}$ from \eqref{eq:commut}, and
$$\left\|u^\e_{\rm free}(t)\right\|_{L^2}=\|f\|_{L^2}\quad ;\quad
\|J^\eps(t)u^\e_{\rm free}\|_{L^2}=\|\nabla f\|_{L^2}.
$$
From the Sobolev inequality (\ref{eq:sobolev}),
$$\left\|u^\e_{\rm free}(t)\right\|_{L^{ s}}\leq \frac{C}{|\cos t|^{\delta
({ s})}} \|f\|_{L^2}^{1-\delta
({ s})}\|\nabla f\|_{L^2}^{\delta
({ s})}$$
for any $s \in [2,\frac{2n}{n-2}[$.
Therefore there exists $C_0$ such that
\begin{equation}\label{eq:estv}
\left\|u^\e_{\rm free}(t)\right\|_{L^{ s}}\leq \frac{C_0}{|\cos t|^{\delta
({ s})}}\, \cdot
\end{equation}
From Prop.~\ref{prop:cauchy}, for fixed $\eps >0$, $u^\eps \in
C(\R,\Sigma)$, and the same obviously holds for $u^\e_{\rm free}$. Therefore,
$w^\eps \in C(\R,\Sigma)$, and there exists $t^\eps>0$ such that
\begin{equation}\label{eq:solong}
\left\|w^\eps(t)\right\|_{L^{s}}\leq \frac{C_0}{|\cos t|^{\delta
({s})}},
\end{equation}
for any $t\in [0,t^\eps]$. So long as (\ref{eq:solong}) holds,
we have
$$\left\|u^\eps(t)\right\|_{L^{ s}}\leq \frac{2C_0}{|\cos t|^{\delta
({s})}},$$
and we can apply Prop.~\ref{prop:estgen}.

Take
$h^\eps=\eps^\gamma \(|x|^{-\gamma}\ast |u^\e|^2\) u^\e_{\rm free}$
and $F^\eps(w^\eps)= \(|x|^{-\gamma}\ast |u^\e|^2\) w^\eps$ and let
$q,k,r,s \in [1,\infty]$ satisfy
the assumptions of Prop.~\ref{prop:estgen}. Now by H\"older's inequality,
\begin{equation*}
\left\|F^\eps(w^\eps)(t)\right\|_{L^{r'}} \leq
\big\| |x|^{-\gamma} \ast |u^\eps(t)|^2\big\|_{L^{\beta}}
\|w^\eps(t)\|_{L^r}
\end{equation*}
with $\beta$ such that $\frac{1}{r'}=\frac{1}{r}+\frac{1}{\beta}$.
By the Hardy--Littlewood--Sobolev inequality and the above estimate,
\begin{equation*}
\begin{split}
\left\|F^\eps(w^\eps)(t)\right\|_{L^{r'}}& \lesssim
\|u^\eps(t)\|^2_{L^{s}}
\|w^\eps(t)\|_{L^r}\\
&\lesssim \frac{(2C_0)^2}{|\cos t|^{2\delta({s})}}
\|w^\eps(t)\|_{L^r}.
\end{split}
\end{equation*}
Note that the second statement of \eqref{eq:holder}(a) ensures that
$s,\beta \in (1,\infty)$ so the
Hardy-Littlewood-Sobolev inequality is applicable here.
Assume (\ref{eq:solong}) holds for $0\leq t\leq T^\e$. If $0 \leq t\leq
T^\e\leq \frac{\pi}{2}-\Lambda \eps$, then $\eps \lesssim \cos
t$, and the
above estimate shows that $F^\eps$ satisfies assumption (\ref{eq:hypF}).

From Corollary~\ref{cor:estgen}, if $\Lambda$ is sufficiently large,
we get for $0\leq t\leq T^\e\leq \frac{\pi}{2}-\Lambda \eps$:
\begin{equation*}\label{eq:wL2}
\|w^\eps\|_{L^\infty(0,T;L^2)}\leq 
C_{\sigma}\eps^{\gamma-1-\frac{1}{\rho}} \left\|(|x|^{-\gamma} \ast
|u^\eps|^2) u^\e_{\rm free}\right\|_{ L^{\rho'}(0,T;L^{\sigma'})}
\end{equation*}
 for any admissible $(\rho,\sigma)$.
Now take $(\rho,\sigma) = (q,r)$ and proceed as above in space, and apply
H\"older inequality in time:
$$\left\| \(|x|^{-\gamma}\ast |u^\e|^2\) u^\e_{\rm free}\right\|_
{ L^{q'}(0,T;L^{r'})}\leq C_{\gamma,n}
\|u^\eps\|^{2}_{L^k(0,T;L^s)}\|u^\e_{\rm free}\|_{ L^q(0,T;L^r)}.
$$
The first term of the right-hand side is estimated through (\ref{eq:estv})
and (\ref{eq:solong}):
$$\left\|u^\eps\right\|^{2}_
{ L^k(0,T;L^{s})}\leq \frac{C}{\left(
\frac{\pi}{2}-T \right)^{2(\delta(s)-1/k)}},
$$
the last term is estimated the same way, for
(\ref{eq:estv}) still holds when replacing $s$ with
$r$:
$$\left\|u^\e_{\rm free}\right\|_
{ L^q(0,T;L^{r})}\leq \frac{C}{\left(
\frac{\pi}{2}-T \right)^{\delta(r)-1/q}}\, .
$$
We infer:
$$\left\| \(|x|^{-\gamma}\ast |u^\e|^2\) u^\e_{\rm free}\right\|_
{ L^{q'}(0,T;L^{r'})}
\leq \frac{C}{\left(
\frac{\pi}{2}-T \right)^{\gamma-1 -\frac{1}{q}}},$$
thus
\begin{equation}\label{eq:reste}
\|w^\eps\|_{L^\infty(0,T;L^2)}\leq 
C  \left(\frac{\eps}{\frac{\pi}{2}-T} \right)^{\gamma -1
-\frac{1}{q}}.
\end{equation}

Now apply the operator $J^\eps$ to (\ref{eq:DuhW}). Since $J^\eps$ and
$\U$ commute, it yields,
\begin{equation*}
J^\eps(t)w^\eps= \U(t)J^\eps(0)r^\eps -i \eps^{\gamma -1}\int_0^t
\U(t-s)J^\eps(s)\left( \(|x|^{-\gamma}\ast |u^\e|^2\) u^\e\right)(s)ds.
\end{equation*}
The action of $J^\e$ on the nonlinear term is described by
\eqref{eq:actionNL}.
In order to apply Prop. \ref{prop:estgen} like before, we take now
$$ h^\eps=\eps^{\gamma-1} \(|x|^{-\gamma}\ast |u^\e|^2\) J^\eps(t)
u^\e_{\rm free}
+ \eps^{\gamma-1}
\ 2 \mbox{Re} \Big(|x|^{-\gamma}\ast  \( \overline{u^\e}
J^\eps(t) u^\e_{\rm free}\) \Big) u^\e
$$
and
\begin{equation} \label{eq:Fgalil}
F^\eps(w^\eps)= \(|x|^{-\gamma}\ast |u^\e|^2\) J^\eps(t)w^\eps
+
\ 2 \mbox{Re} \Big(|x|^{-\gamma}\ast \( \overline{u^\e} J^\eps(t) w^\e \)\Big) u^\e \, .
\end{equation}
The first term on the r.h.s of \eqref{eq:Fgalil} leads to an
equation which is very similar to (\ref{eq:DuhW}), with $w^\eps$
replaced by $J^\eps w^\eps$ 
and is treated by the same computations as above.
For the second term, we estimate by H\"older, by
the Hardy-Littlewood-Sobolev
inequality and then again by H\"older:
\begin{equation*}
\begin{split}
\left\|
 2 \mbox{Re} \Big(|x|^{-\gamma}\ast \( \overline{u^\e} J^\eps(t) w^\e \)\Big) u^\e
\right\|_{L^{r'}} & \lesssim
\| |x|^{-\gamma}\ast \(
\overline{u^\e} J^\eps(t) w^\e \)  \|_{L^{\beta_1}}\|u^\eps(t)\|_{L^{s}} \\
& \lesssim
\|u^\eps(t)\|_{L^{s}}\| J^\eps(t) w^\eps(t)\|_{L^r}\|u^\eps(t)\|_{L^{s}}
\end{split}
\end{equation*}
with $r,s$ as stated in \eqref{eq:holder} and
$\frac{1}{r'} = \frac{1}{\beta_1}+\frac{1}{s}$ .
Here the condition to use the Hardy-Littlewood-Sobolev
inequality is $\gamma > \delta(r)+\delta(s)$, which is always satisfied by
 \eqref{eq:holder}.
By applying now \eqref{eq:solong} we continue to estimate
$$
\leq \frac{(2C_0)^2}{\left(|\cos t| \right)^{2\delta({s})}}
\| J^\eps(t) w^\eps(t)\|_{L^r}\, .
$$
Then we apply, like before,
Prop. \ref{prop:estgen} and estimate the term for $h^\e$ as above,
to obtain:
\begin{equation}\label{eq:reste2}
\|J^\eps w^\eps\|_{L^\infty(0,T;L^2)}\leq C\|\nabla r^\eps\|_{L^2}+ C
\left(\frac{\eps}{\frac{\pi}{2}-T} \right)^{\gamma -1
-\frac{1}{ q}}.
\end{equation}
Combining (\ref{eq:reste}) and (\ref{eq:reste2}) yields, along with
(\ref{eq:sobolev}),
$$ \forall t\in [0,T],\ \  \left\|w^\eps(t)\right\|_{L^{s}}
\leq \frac{C}{|\cos t|^{\delta
({s})}}\left(\|r^\eps\|_{H^1} +
\left(\frac{\eps}{\frac{\pi}{2}-t} \right)^{\gamma -1
-\frac{1}{q}} \right).$$
Therefore, choosing $\eps$ sufficiently small and $\Lambda$
sufficiently large, we deduce that we can take
$T=\frac{\pi}{2}-\Lambda \eps$. With the result of Lemma \ref{lem:DA.1}
on the limit of $u^\e_{\rm free}$,
this yields Prop.~\ref{prop:scattering},
away from the focus,
for
$A^\eps=Id$ and $J^\eps$. The case $A^\eps=H^\eps$ on this
time interval is now
straightforward.

\medskip

The remaining parts of the proof for Prop.~\ref{prop:scattering}
are done as in \cite{CaIHP} with the method changes as in the
part shown above.
It remains to show that the  approximations in the two
different regimes match  at
$t_\ast=\frac{\pi}{2}-\Lambda \eps$,
and that the influence of the harmonic potential is small near the focus
so that the propagation there is given by
\begin{equation}\label{eq:r3gscal}
 v^\eps(t,x)  = \frac{1}{\e^{n/2}}\ \psi \( \frac{t-\frac{\pi}{ 2}}{\eps} ,
\frac x \eps \) ,
\end{equation}
where $\psi$ is the solution of \eqref{eq:r3gen} subject to the
following initial condition at $t=-\infty$
$$
 {\tt U}(-t)\psi(t)\big|_{t=-\infty} =
 e^{i\frac{n\pi}{4}}\widehat{f}.
$$
This solution exists according to Proposition \ref{prop:scattr3}.

Then the following asymptotic is proven:
\begin{equation*}
\limsup_{\eps \to 0} \sup_{\frac \pi 2 - \Lambda \e \leq t \leq \frac
\pi 2 + \Lambda \e}
\Big\| A^\e(t)\(u^\e(t,x) -  v^\eps(t,x)  \)
\Big\|_{L^2_x}\Tend{\Lambda }{+\infty}  0\, ,
\end{equation*}
with $A^\e(t)$ being one of the operators  $Id$, $J^\e$ or $H^\e$.
Since these parts are quite similar to the treatment in \cite{CaIHP},
we do not repeat them.

After the crossing of the first focus, the solution is again propagated
linearly and at subsequent focusing points this process is iterated.
 \qed


\section{Formal Computations and Discussions}
\label{sec:formal}
\subsection{The case $\alpha = \gamma = 1$ (in 3-d: Schr\"odinger-Poisson)}
We saw in Section~\ref{sec:scattering} that when $\alpha =\gamma>1$,
the nonlinear term in \eqref{eq:r3.3} has a leading order influence
near the focuses, and only in these regions. On the other hand, if
$\alpha=1$ and $\gamma <1$, Section~\ref{sec:nlwkb} shows that the
Hartree term cannot be neglected away from the focuses. These two
cases suggest that when $\alpha =\gamma =1$, the nonlinear influence
is everywhere relevant. The aim of this final section is to give
convincing arguments that this is the case.

\smallskip

For the influence near the focuses, we need the scattering theory
for  \eqref{eq:r3gen} at $\gamma = 1$. In this
long range scattering case,  modified scattering operators
are needed instead of the ones described in Prop. \ref{prop:scattr3}.
Hayashi and Naumkin \cite{HN98} obtained an asymptotic completeness result
for $n\geq2$ with smoothness assumptions which are applicable to
our situation.
On the other hand, they could not obtain wave operators.
Ginibre and Velo \cite{GVII,GVIII} obtained modified wave operators
for  \eqref{eq:r3gen} with $\gamma = 1$ using Gevrey spaces by a
technically involved method. A drawback of both these results is
that they include a loss in regularity.

To show how the long range
scattering theory fits into our framework we report (a particular case
of) the result of
Hayashi and Naumkin \cite{HN98}.
\begin{prop}[\cite{HN98}]\label{prop:longrange}
Assume $n=3$, $\varphi \in \Sigma$,
 and  $\delta = \|\varphi\|_\Sigma$ is sufficiently small.
Let $\psi\in C(\R,\Sigma)$ be the solution
of \eqref{eq:r3gen} with $\psi_{t=0}=\varphi$. Then
there exists a unique function $\psi_+ \in  H^{\sigma,0}\cap H^{0,\sigma}$,
$\frac 1 2 <\sigma<1$, such that
$$
\left\|
\psi(t) - \exp\( i
\(|x|^{-1}\ast|\widehat  \psi_+|^2\)
\(\frac x t\)\log |t| \) {\tt U}(t) \psi_+
\right\|_{L^2} \Tend{t }{+\infty} 0\, ,
$$
where  $H^{\alpha,\beta} = \{\phi \in  \mathcal S' \ \big| \ \|
(1+|x|^2)^{\beta/2}
(1+\Delta)^{\alpha/2}\phi\|_{L^2}<\infty\} $.
\end{prop}
To summarize very roughly, the results in \cite{GVII,GVIII} consist in
showing that given some $\psi_+$ (or $\psi_-$ for an asymptotic
behavior for $t\to -\infty$), one can find $\psi$ solving
\eqref{eq:r3gen} such that the above asymptotics holds.

Analogoulsy to the treatment of long-range scattering in \cite{CaCMP},
one can now define
$ g^\e(t,x) := \(|x|^{-1}\ast| f|^2\)(x) \log(
\frac{\cos t}{\e})$
(compare with \eqref{eq:g})
and add the phase $g^\e\bigl|_{t=0}$  to the initial
data in \eqref{eq:r3.3}. This yields:
\begin{equation*}
u^\e\big|_{t=0} = f(x) e^{-i\(|x|^{-1}\ast| f|^2\)(x) \log
\e}\, .
\end{equation*}

Using the modified scattering operators from the
results of \cite{GVII,GVIII} we get, at least formally,
for $0\leq t <\pi/2$,
\begin{equation*}\label{eq:formpro}
 u^\e(t ,x)\sim
\frac{1}{(\cos t)^{3/2}}
 f\(\frac{x}{\cos t}\) e^{-i\frac{x^2}{2\e}\tan t +ig^\e\(
t,\frac{x}{\cos t}\) }
\quad \text{as }\e \to 0\,.
\end{equation*}
This asymptotic also stems from the same computations as those
performed in Section~\ref{sec:osc}.
Notice that the matching for $|t-\frac{\pi}{2}|=\O(\e)$ is similar to the one
in \cite{CaIHP}, except that we now have to take the presence of
$g^\e$ into account. This is where changing the integration from $0$
to $t$ in \eqref{eq:g} into the above definition of $g^\e$ makes the
matching possible. Indeed,  for $|t-\frac{\pi}{2}|=\O(\e)$, we compare
$u^\e$ with the
function $v^\e$ given by \eqref{eq:r3gscal}, where $\psi$ is now the
solution given by the long range wave operators constructed in
\cite{GVII,GVIII}. To make this statement more precise and the link
between \eqref{eq:g} and the definition of $g^\e$ more explicit,
notice that we have, as $t\to \frac{\pi}{2}$:
\begin{align*}
g^\e\(t,\frac{x}{\cos t}\)&\sim \(|x|^{-1}\ast| f|^2\)\(\frac{x}{\frac{\pi}{2}
- t}\) \log\( \frac{\frac{\pi}{2}-t}{\e}\) \quad \text{(phase shift
for }v^\e )\\
&\sim -\(|x|^{-1}\ast|
f|^2\)\(\frac{x}{\cos
t}\)\int_{\arccos \e}^t \frac{d\tau}{\cos \tau}\quad \text{(compare
with \eqref{eq:g})}\, .
\end{align*}

  The effects of the nonlinearity show up in  $g^\e$.
Using the scaling \eqref{eq:r3gscal} we can then (formally)
continue with Prop.~\ref{prop:longrange}: for $\pi/2< t <3\pi/2$,
\begin{equation*}
 u^\e(t ,x)\sim
\frac{e^{-i\frac{3 \pi}{2}}}{|\cos t|^{3/2}}
\(\F\circ \widetilde S\circ\F^{-1}\)  f\(\frac{x}{\cos t}\)
e^{-i\frac{x^2}{2\e}\tan t  +ih^\e\(
t,\frac{x}{\cos t}\) }\quad \text{as }\e \to 0\, ,
\end{equation*}
where $\widetilde S$ is the map $\widetilde S:\psi_-\mapsto \psi_+ $,
where $\psi_-$ is the asymptotic state of the result of \cite{GVII},
which yields some solution $\psi$ to \eqref{eq:r3gen}, and $\psi_+$
is provided by Prop. \ref{prop:longrange}.
$h^\e$ is given by
$$ h^\e(t,x) := -\(|x|^{-1}\ast| \F\circ \widetilde S\circ\F^{-1} f|^2\)(x)
\log\( \frac{|\cos t|}{\e}\)\, .$$
The action of $\F\circ \widetilde S\circ\F^{-1}$ on $f$ accounts for
nonlinear effects taking place at the focus,
and $h^\e$ for nonlinear effects after the focus.
So the influence of the nonlinearity will be relevant at all times.

The impossibility to define a scattering operator for this
case  is one of the reasons why this argument is only formal.

\begin{rema*}
A rigorous result could  be obtained with the same approach as in
\cite{Ca2}. It would consist in studying the system of \emph{linear}
equations with a \emph{nonlinear coupling},
\begin{equation*}\left\{
\begin{aligned}
i\e\d_t {\bf u}^\e +\frac{1}{2}\e^2\Delta {\bf u}^\e &=
\frac{|x|^2}{2}{\bf u}^\e \, ,\\
i\e\d_t u^\e +\frac{1}{2}\e^2\Delta u^\e &=
\frac{|x|^2}{2}u^\e +\e \( |x|^{-1}\ast |{\bf u}^\e|^2\)u^\e\, .
\end{aligned}\right.
\end{equation*}
The first equation is solved explicitly thanks to Mehler's formula,
and the second one is a linear Schr\"odinger equation with a harmonic
potential and a time-dependent perturbation. With the oscillatory
integral used in Section~\ref{sec:nlwkb}, and adapting the results of
\cite{DG}, one could prove similar asymptotics to those stated above.
\end{rema*}

\subsection{The case of an additional local strong nonlinearity}
We now consider equation \eqref{eq:r3.3} with an additional nonlinear
term that is a multiplication operator with a power of the density
$|\u^\e|^2$.

Such equations arise in the modeling of effective one particle
Schr\"odinger equations where ``exchange terms" like in the Hartree-Fock
equation are simplified to functionals of the local densities, i.e.
time dependent density functional theory, with the Schr\"odinger-
Poisson-X$\alpha$ equation as the simplest of such models (see
\cite{MauserAML}
and \cite{BMS1} for a heuristic derivation and numerical simulations).
Note that the additional ``local" term has the opposite sign than the
Hartree term (corresponding to the physical fact that the
``exchange-correlation hole'' weakens the direct Coulomb interaction).

We will hence consider the following class of semi-classical Hartree equations
\begin{equation}\label{eq:r3.3-Xa}
i\e \d_t  u^\e +\frac{1}{2}\e^2\Delta u^\e = \frac{|x|^2}{2}u^\e +
\e^\alpha \(|x|^{-\gamma}\ast |u^\e|^2\)u^\e - \e^\beta |u^\e|^{2\sigma} u^\e ,
\end{equation}
with  $\alpha \geq 1$, $\beta \geq 1$, $\gamma>0$ for $x\in\R^n$,
and with a $\sigma$ that is sub-critical with respect to finite time blow-up,
i.e. $\frac{2}{n}>\sigma>0$.

We can now discern
the influence of the two nonlinear terms in the classical
limit in terms of:
\begin{itemize}
\item The size of the scaling exponents $\alpha$, $\beta$ with respect
to the critical value .
\item The relation between the scaling and the ``strength'' of the
nonlinearities determined respectively by $\gamma$ and $\sigma$.
\end{itemize}
If we take $\alpha>1$ and $\beta>1$, by
\cite{CaIHP} and Section \ref{sec:scattering}
we find that the classical limit is given by the linear propagation as long as
no focusing occurs. At the focus, the relevant discrimination is
$\sigma=\beta/n$ or $<\beta/n$ for the power nonlinearity and
$\gamma=\alpha$ or $<\alpha$ for the Hartree term.
If $\sigma=\beta/n$ and $\gamma<\alpha$, the crossing of the focus
will be described by the scattering operator for NLS (when it is defined),
 if on the other hand
$\sigma<\beta/n$ and $\gamma=\alpha$
(and the assumptions of Prop.~\ref{prop:scattr3} are satisfied),
 focus crossing will be determined by
the scattering operator of Prop.~\ref{prop:scattr3}.
If both nonlinearities are at the critical strength
($\sigma=\beta/n$ and $\gamma=\alpha$),
 then both will have an influence in crossing
the caustic. If, on the other hand, both $\sigma<\beta/n$ and $\gamma<\alpha$,
the nonlinear influence will be negligible everywhere.

If at least one of the scaling exponents $\alpha$ and $\beta $ is equal to
1 and, at the same time, both $\sigma<\beta/n$ and $\gamma<\alpha$,
the corresponding nonlinear term will be relevant in the WKB propagation
before the focusing. At the focus, the nonlinear terms will not be
relevant and the crossing of the focus will be as in Prop. \ref{prop:nlwkb}.
If  $\sigma=\beta/n$ and $\gamma=\alpha$ then there will be a nonlinear
influence everywhere and long range scattering for NLS and/or
Hartree has to be taken into account.

The influence of the nonlinear action for the single power NLS
and the Hartree equation is summed up in two tables, for Hartree
the table is given in the introduction, for single power nonlinear Schr\"odinger equation,
it is stated in \cite{Ca2}.
The behavior of \eqref{eq:r3.3-Xa} can be described by independently
superposing these two tables. The following table
is an extract from that superposition:
\begin{center}
\begin{tabular}[c]{l|c|c}
&&\\
 &$\alpha >\gamma$ {\bf and} $\beta > \sigma n$
 & $\alpha =\gamma$ {\bf or} $\beta = \sigma n$\\
\hline
&&\\
$\alpha >1$ {\bf and} $\beta> 1$ &Linear WKB, & Linear WKB,\\
&linear focus   & nonlinear focus \\
&&\\
\hline
&&\\
$\alpha =1$ {\bf or} $\beta=1$&Nonlinear WKB, &Nonlinear WKB,\\
&linear focus   & nonlinear focus \\
\end{tabular}
\end{center}
``Nonlinear WKB'' respectively ``nonlinear focus''
here stands for an influence from at least
one of the nonlinear terms  away from the focus or close to the focus.

\subsection{Wigner measures}
We already mentioned in the introduction the work of Zhang, Zheng and Mauser
\cite{ZZM} where the (semi)classical limit of the Schr\"odinger-Poisson
equation with no smallness assumption (on the initial data or the nonlinearity)
is studied by means of Wigner measures.
Wigner measures have proven to be efficient tools for linear
semi-classical problems and for homogenization limits; see \cite{MauserForges}
for an overview  on Wigner measure limits of Hartree equations.
Wigner measures have the merit that in phase space the caustics of physical
space are somewhat unfolded and that generally, results globally in time
are possible. \\
In \cite{CaWigner}, the  Wigner measure of the nonlinear Schr\"odinger equation
with power-like nonlinearity studied in \cite{Ca2} is investigated.
It is shown that the Wigner measure
leads to an ill-posed problem whenever nonlinear effects at the focal
points come into play. In other words, the Wigner measure can only be
valid as long as no caustic appears. We  briefly discuss the
Wigner measures of  \eqref{eq:r3.3} in view of these results.

The Wigner measure of the family $(u^\e(t))_{0<\e\leq1}$, which is bounded
in $L^2$, is the weak limit under $\e \to 0$ (up to an extraction) of its Wigner transform,
\begin{equation*}
W^\e(u^\e)(t,x,\xi)=
\frac{1}{(2\pi)^n}
\int u^\e\(t,x-\frac{v\e}{2}\)   \overline{u^\e}\(t,x+\frac{v\e}{2}\)
e^{i\xi \cdot v} dv.
\end{equation*}
This limit is a positive radon measure $\mu$ and is in general not a unique
limit.

-- {\em linear case:} Case $\alpha >\gamma$, $\alpha > 1$:\\
By the result of Prop. \ref{prop:small} and the asymptotics of $u_{\rm free}$
in Lemma \ref{lem:DA.1},  the Wigner measure $\mu^-$
for $t<\pi/2$ of the family $(u^\e(t))_{0<\e\leq1}$ is
$$
\mu^-(t,x,\xi)=\frac{1}{|\cos t|^{n}} \left|
f\(\frac{x}{\cos t}\)\right|^2 dx \otimes \delta_{\xi=x \tan t}.
$$
For $\pi/2 < t < \pi$, the Wigner measure  of $(u^\e(t))_{0<\e\leq1}$
(denoted by $\mu^+$) is the same:
$\mu^+(t,x,\xi) = \mu^-(t,x,\xi)$. At $t=\pi/2$,
the limits from above and below are:\\
$\lim_{t\to\pi/2^-}\mu^-(t,x,\xi) =\lim_{t\to\pi/2^+}\mu^+(t,x,\xi) =
|f(\xi)|^2 d\xi \otimes \delta(x)$.

-- {\em nonlinear WKB, linear focus:} Case $\gamma < \alpha = 1$. \\
The asymptotics of $u^\e$ are stated in Prop. \ref{prop:nlwkb}.
The additional phase term $g$
is of order $1$ and does not change the Wigner measure of
$(u^\e(t))_{0<\e\leq1}$, so in this case  $\mu^-$ and $\mu^+$ are the
same as in the previous case: the Wigner measure does not "see" the nonlinear effect $g$.

-- {\em linear WKB, nonlinear focus:} Case $\gamma = \alpha > 1$. \\
The asymptotics of Prop. \ref{prop:scattering} involve, for $t\geq\pi/2$,
the scattering operator associated with the unscaled equation
\eqref{eq:r3gen}. For  $t<\pi/2$, the Wigner measure of
$(u^\e(t))_{0<\e\leq1}$ is still the same as above, but for
$\pi/2<t<\pi$, we have
$$
\mu^+(t,x,\xi)=\frac{1}{|\cos t|^{n}} \left|
\F\circ S\circ\F^{-1}
f\(\frac{x}{\cos t}\)\right|^2 dx \otimes \delta_{\xi=x \tan t},
$$
where $S$ is the scattering operator for \eqref{eq:r3gen} and $\F$
the Fourier transform.

-- {\em nonlinear WKB, nonlinear focus:} Case $\gamma = \alpha = 1$. \\
The asymptotics for this case of (the formal computation)
Prop. \ref{prop:longrange} include an additional phase term
which is of order $\log \e$ and a modification of the initial data of the same
order of magnitude. Both do not alter the Wigner measure,
since they are dominated
by the scaling of the Wigner transform, and thus the Wigner measure
 is the same as in the previous case.\\
For the last two cases, the limits at $t=\pi/2$ are
\begin{equation}\label{eq:wigdis}
\begin{split}
&\lim_{t\to\pi/2^-}\mu^-(t,x,\xi) = |f(\xi)|^2 d\xi \otimes \delta(x),
\\
&\lim_{t\to\pi/2^+}\mu^+(t,x,\xi) = |\F\circ S\circ\F^{-1}f(\xi)|^2 d\xi \otimes \delta(x).
\end{split}
\end{equation}
The idea of \cite{CaWigner} is to find now two profiles $f_1$ and $f_2$
for which $|f_1|^2 \equiv |f_2|^2$, but at the same time
$|\F\circ S\circ\F^{-1}f_1|^2\not \equiv|\F\circ S\circ\F^{-1}f_2|^2$ .
Then the Wigner measures of the corresponding families
$(u_j^\e(t))_{0<\e\leq1}$, $j=1,2$, will be equal up to the focus,
but different after the focus, i.e.
$\mu_1^-=\mu_2^-$ but $\mu_1^+ \neq \mu_2^+$.
So after the caustic point the Wigner measure will not be unique
anymore in the case where the
nonlinearity is relevant at the focus.
These profiles were constructed using an expansion of $S$ around the origin.
Since our problem is very similar to the one studied there, we expect a
similar result to hold for equation \eqref{eq:r3.3},
i.e. we expect the Wigner measure to lead to
an ill-posed problem if there is a nonlinear influence at the caustic.

In view of the result of \cite{ZZM},
note that the non-uniqueness of the weak solutions for Vlasov-Poisson with
measures as initial data and the non-uniqueness of the Wigner measure of
a given $\eps$-dependent family of solutions coincide, such that there is
no contradiction with the global and unique semi-classical limits
of the
Hartree type equations obtained here.

\bibliographystyle{amsplain}
\bibliography{carles}

\end{document}